\documentclass[onecolumn,prl]{revtex4}

\usepackage{graphicx}
\usepackage{type1cm}
\usepackage{subfigure}
\usepackage{amssymb,latexsym,amsmath,epsfig,subfig,epic,color,amscd,dsfont,lettrine}
\usepackage{oldgerm,mathrsfs,euscript,setspace,bbm}
\usepackage{soul}

{

}

\newcommand{\R}{\mathbb{R}}

\newcommand{\Ed}{\mathcal{E}}

	\usepackage{makeidx}
	\usepackage[ansinew]{inputenc}
	\usepackage[usenames,dvipsnames]{pstricks}
	\usepackage{pst-grad} 
	\usepackage{pst-plot} 

\makeindex

\begin{document}

\title{Limitations and tradeoffs in synchronization of large-scale networks with uncertain links}

\author{Amit Diwadkar and Umesh Vaidya}
\affiliation{Electrical and Computer Engineering, Iowa State University \\ Coover Hall, Ames, IA, U.S.A. 50011 \\
Correspondence to diwadkar@iastate.edu, ugvaidya@iastate.edu}

\begin{abstract}
The synchronization of nonlinear systems connected over large-scale networks has gained popularity in a variety of applications, such as power grids, sensor networks, and biology. Stochastic uncertainty in the interconnections is a ubiquitous phenomenon observed in these physical and biological networks. We provide a size-independent network sufficient condition for the synchronization of scalar nonlinear systems with stochastic linear interactions over large-scale networks. This sufficient condition, expressed in terms of nonlinear dynamics, the Laplacian eigenvalues of the nominal interconnections, and the variance and location of the stochastic uncertainty, allows us to define a synchronization margin. We provide an analytical characterization of important trade-offs between the internal nonlinear dynamics, network topology, and uncertainty in synchronization. For nearest neighbour networks, the existence of an optimal number of neighbours with a maximum synchronization margin is demonstrated. An analytical formula for the optimal gain that produces the maximum synchronization margin allows us to compare the synchronization properties of various complex network topologies.
\end{abstract}

\maketitle

\lettrine[lines=2,findent=1pt]{S}{}ynchronization in large-scale network systems is a fascinating problem that has attracted the attention of researchers in a variety of scientific and engineering disciplines. It is a ubiquitous phenomenon in many engineering and naturally occurring systems, with examples including generators for electric power grids, communication networks, sensor networks, circadian clocks, neural networks in the visual cortex, biological applications, and the synchronization of fireflies \cite{Strogatz2, kuramoto_bonilla1,genetic_clock_hasty,PNAS_Dorfler}. The synchronization of systems over a network is becoming increasingly important in power system dynamics. Simplified power system models demonstrating synchronization are being studied to gain insight into the effect of network topology on the synchronization properties of dynamic power networks \cite{Chaos_Rohden}. The effects of network topology and size on the synchronization ability of complex networks is an important area of research \cite{Chaos_Stout}. Complex networks with certain desirable properties, such as a small average path between nodes, low clustering ability, and the existence of hub nodes, among others, have been extensively studied over the past decade \cite{Erdos_Renyi,Nature_watts,PNAS_Amaral,Science_Barabasi,Pecora_Nature,Arndt_Nature}.

It is impossible to do justice to the long list of literature that exists in the area of synchronization of dynamical systems. In the following discussion, we list a few references that are particularly relevant to the results presented in this paper. In \cite{pecora_carroll_master}, the master stability function was introduced to study the local synchronization of chaotic oscillator systems. Interesting computational observations were made that indicated the importance of the smallest and largest eigenvalues of the graph Laplacian. The master stability function was also used to study synchronization over Small-World networks and provide bounds on the coupling gains to guarantee the stability of the synchronous state in \cite{pecora_barahona}. Bounds were provided on the coupling gains to guarantee the stability of the synchronous state in \cite{rangarajan_ding}. The impact of network interconnections on the stability of the synchronous state of a network system was also studied in \cite{belykh_connection}. These results derived a condition for global synchronization based on the coupling weights and eventual dissipativity of the chaotic system using Lyapunov function methods and a bound on path lengths in the connection graph.  In this paper, as in the papers listed above, we provide an analytical characterization of the importance of the smallest and largest positive eigenvalue of the coupling Laplacian. However, in contrast to the above references, we provide conditions for the global synchronization in the presence of stochastic link uncertainty. Understanding the role of spatial perturbation in the nearest neighbour network to force a transition from one synchronized state to another is important for molecular conformation \cite{PNAS_Mezic}. Other aspects of network synchronization that are gaining attention are the effects of network topology and interconnection weights on the robustness of the synchronization properties \cite{PNAS_Nishikawa}. In this paper, we provide a systematic approach for understanding the effects of stochastic spatial uncertainties, network topology, and coupling weights on network synchronization.

Uncertainty is ubiquitous in many of these large-scale network systems. Hence, the problem of synchronization in the presence of uncertainty is important for the design of robust network systems. The study of uncertainty in network systems can be motivated in various ways. For example, in electric power networks, uncertain parameters or the outage of transmission lines are possible sources of uncertainty. Similarly, a malicious attack on network links can be modelled as uncertainty. Synchronization with limited information or intermittent communication between individual agents, e.g., a network of neurons, can also be modelled using time-varying uncertainty. In this paper, we address the problem of robust synchronization in large-scale networks with stochastic uncertain links. Existing literature on this problem has focused on the use of Lyapunov function-based techniques to provide conditions for robust synchronization \cite{synchronization_scalefree}.

Both the master stability function and Lyapunov exponents have been used to study the variation of the synchronous state's stability, given local stability results with stochastic interactions \cite{porfiri_master,amit_synchronization}. The problem of synchronization in the presence of simple on-off or blinking interaction uncertainty was studied in \cite{belykh_blinking,Belykh_SIAM_I,Belykh_SIAM_II} using connection graph stability ideas \cite{belykh_connection}. The local synchronization of coupled maps was studied in \cite{Jost_SIAM,Jost_Atay}, which also provides a measure for local synchronization. Synchronization over balanced neuron networks with random synaptic interconnections has also been studied \cite{PRE_Molino}. Researchers have studied the emergence of robust synchronized activity in networks with random interconnection weights \cite{PRE_Sinha}. The robustness of synchronization to small perturbations in system dynamics and noise has been studied \cite{PRE_Kocarev}, while the robustness to parameter variations was also studied in the context of neuronal behaviour \cite{PRE_Wang}. In this paper, we consider a more general model for stochastic link uncertainty than the simple blinking model and develop mathematically rigorous measures to capture the degree of synchronization.

We consider a network of systems where the nodes in the network are dynamic agents with scalar nonlinear dynamics. These agents are assumed to interact linearly with other agents or nodes through the network Laplacian. The interactions between the network nodes are assumed to be stochastic. This research builds on our past work, where we developed an analytical framework using system theoretic tools to understand the fundamental limitations of the stabilization and estimation of nonlinear systems with uncertain channels \cite{amit_erasure_observation_journal,amit_ltv_journal,Vaidya_erasure_SCL,Vaidya_erasure_stablization}. There are two main objectives for this research, which also constitute the main contributions of this paper. The first objective is to provide a scalable computational condition for the synchronization of large-scale network systems. We exploit the identical nature of the network agent dynamics to provide a sufficient condition for synchronization, which involves verifying a {\it scalar} inequality. This makes our synchronization condition independent of network size and hence computationally attractive for large-scale network systems. The second objective and contribution of this paper is to understand the {\it interplay} between three network characteristics: (1) internal agent dynamics, (2) network topology captured by the nominal graph Laplacian, and (3) uncertainty statistics in the network synchronization. We use tools from robust control theory to provide an analytical expression for the synchronization margin that involves all three network parameters and increases the understanding of the {\it trade-offs} between these characteristics and network synchronization. This analytical relationship provides useful insight and can compare the robustness properties for nearest neighbour networks with varying numbers of neighbours. In particular, we show that there exists an optimal number of neighbours in a nearest neighbour network that produces a maximum synchronization margin. If the number of neighbours is above or below this optimal value, then the margin for synchronization decreases.

We use an analytical expression for the optimal gain and synchronization margin to compare the synchronization properties of Small-World and Erdos-Renyi network topologies.
\section{Results}
\subsection{Synchronization in Dynamic Networks with Uncertain Links}\label{section_main}
We consider the problem of synchronization in large-scale nonlinear network systems with the following scalar dynamics of the individual subsystems:
\begin{align} \label{component_dynamics}
x^k _{t+1} = a x^k_t - {\phi}  (x^k_t) + v^k_t \qquad ~ ~k=1,\ldots, N,
\end{align}
where $x^k \in \mathbb{R}$ are the states of the $k^{th}$ subsystem and $a > 0$ and $v^k \in \mathbb{R}$ is an independent, identically distributed (i.i.d.) additive noise process with zero mean (i.e., $E[v^k_t] = 0$) and variance $E[(v^k_t)^2] = \omega^2$. The subscript $t$ used in Eq. \eqref{component_dynamics} denotes the index of the discrete time-step throughout the paper. The function $\phi \colon \mathbb{R} \to \mathbb{R} $ is a monotonic, globally Lipschitz function with $\phi(0) = 0$ and Lipschitz constant $\frac{2}{\delta}$ for $\delta > 0$.

The individual subsystem model is general enough to include systems with steady-state dynamics that could be stable, oscillatory, or chaotic in nature. We assume the individual subsystems are linearly coupled over an undirected network given by a graph $G=(V,\Ed)$ with node set $V$, edge set $\Ed$, and edge weights $\mu_{ij} \in \mathbb{R}^+$ for $i,j \in V$ and $e_{ij}\in \Ed$. Let $\Ed_U \subseteq \Ed$ be a set of uncertain edges and $\Ed_D = \Ed\setminus \Ed_U$. The weights for $e_{ij} \in \Ed_U$ are random variables: $\zeta_{ij} = \mu_{ij}+\xi_{ij}$, where $\mu_{ij}$ models the nominal edge weight and $\xi_{ij}$ models the zero-mean uncertainty ($E[\xi_{ij}] = 0$) with variance $E[\xi_{ij}^2] = E\left[(\zeta_{ij}-\mu_{ij})^2\right] = \sigma_{ij}^2$. Because the network is undirected, the Laplacian for the network graph is symmetric. We denote the nominal graph Laplacian by ${\cal L} := \left[ l(ij)\right] \in \R^{N \times N}, \; e_{ij}\in \Ed$, where $l(ij) = -\mu_{ij}$, if $i\neq j$, and, $e_{ij} \in \Ed$, $l(ij) = \sum_{e_{ij}\in \Ed}\mu_{ij}$, if $i=j$. We denote the zero-mean uncertain graph Laplacian by ${\cal L}_R := \left[l_R(ij)\right] \in \R^{N \times N}, \; e_{ij}\in \Ed_U$, where $l_R(ij) =  -\xi_{ij}$, if $i\neq j$, and, $e_{ij} \in \Ed_U$, $l_R(ij) = \sum_{e_{ij}\in \Ed_U}\xi_{ij}$, if $i=j$. The nominal graph Laplacian ${\cal L}$ is a sum of the graph Laplacian for the purely deterministic graph $(V,\Ed_D)$, and of the mean Laplacian for the purely uncertain graph $(V,\Ed_U)$. Hence, ${\cal L}$ may be written as ${\cal L} = {\cal L}_D + {\cal L}_U$, where ${\cal L}_D$, is the Laplacian for the graph over $V$ with edge set $\Ed_D$. ${\cal L}_U$ is the mean Laplacian for the graph over $V$ with edge set $\Ed_U$. Define $\tilde x_t  = [ (x^1_t)\; \cdots\; (x^N_t)  ]^{\top} \in \mathbb{R}^N$ and $\tilde {\phi}(\tilde x_t) = [ ({\phi}^1_t)(x^1_t)\;\cdots\; ({\phi}^N_t)(x^N_t) ]^{\top} \in \mathbb{R}^N$, where $A^{\top}$ denotes the transpose of matrix $A$. In compact form, the network dynamics are written as
\begin{align}
\tilde{x}_{t+1} &= \left(aI_N - g({\cal L}+{\cal L}_R)\right) \tilde{x}_t - \tilde{\phi} \left (\tilde x_t\right) + \tilde{v}_t, \label{coupled_dynamics}
\end{align}
where $g > 0$ is the coupling gain and $I_N$ is the $N\times N$ identity matrix.  Our objective is to understand the interplay of the following network characteristics: the internal dynamics of the network components, the network topology, the uncertainty statistics, and the coupling gain for network synchronization. Given the stochastic nature of network systems, we propose the following definition of mean square synchronization \cite{Hasminskii_book}.

\noindent{\bf Mean Square Synchronization:} Define $\Xi := \{\xi_{ij}\;|\;e_{ij}\in \Ed_U\}$ and $E_{\Xi}[\cdot]$ as the expectation with respect to uncertainties in the set $\Xi$. The network system \eqref{coupled_dynamics} is said to be mean square synchronizing (MSS) if there exist positive constants $\beta < 1$, $\bar K <\infty$, and $L <\infty$, such that
\begin{eqnarray}\label{MSS_def}
& E_{\Xi} \parallel x^k_{t} - x^j_{t}\parallel^2 \leq \bar K {\beta}^t \parallel x^k_{0} - x^j_{0}\ \parallel^2 + L\omega^2,
\label{mse_sync_eqn}
\end{eqnarray}
 $\forall k,j \in [1,N]$, where $\bar K$ is a function of $\parallel x_0^i-x_0^j\parallel^2$ for $i,j \in [1,N]$ and $\bar K(0)=K$ is a constant. In the absence of additive noise $\tilde v_t$ in system Eq. \eqref{coupled_dynamics}, the term  $L \omega^2$ in Eq. \eqref{mse_sync_eqn} vanishes and the system is mean square exponential (MSE) synchronizing \cite{Wang_IEEE_TNN}.
We introduce the notion of the {\it coefficient of dispersion} to capture the statistics of uncertainty.

\noindent{\bf Coefficient of Dispersion:} Let $\zeta \in \mathbb{R}$ be a random variable with mean $\mu > 0$ and variance $\sigma^2 > 0$. The coefficient of dispersion (CoD) $\gamma$ is defined as $\gamma := \frac{\sigma^2}{\mu}$. For all edges $(i,j)$ in the network, the mean weights assigned are positive, i.e., $\mu_{ij} > 0$ for all $(i,j)$. Furthermore, the CoD for each link is given by $\gamma_{ij} = \frac{\sigma_{ij}^2}{\mu_{ij}}$ and $\bar{\gamma} = \underset{\xi_{ij}}{\max}{\gamma_{ij}}$.

Because the subsystems are identical, the synchronization manifold is spanned by the vector $\mathds{1}=[1,\ldots,1]^{\top}$. The dynamics on the synchronization manifold are decoupled from the dynamics off the manifold and are essentially described by the dynamics of the individual system, which could be stable, oscillatory, or complex in nature. We apply a change of coordinates to decompose the system dynamics on and off the synchronization manifold.
Let ${\cal L} = V\Lambda V^{\top}$, where $V$ is an orthonormal set of vectors given by $V = \left[\frac{\mathds{1}}{\sqrt{N}}\; U\right]$, in which $U$ is a set of $N-1$ orthonormal vectors that are orthonormal to $\mathds{1}$. Furthermore, we have $\Lambda = \text{diag}\{\lambda_1,\cdots,\lambda_N\}$, where $0=\lambda_1 < \lambda_2 \leq \cdots \leq \lambda_N$ are the eigenvalues of ${\cal L}$. Let $\tilde z_t = V^{\top}\tilde x_t$ and $\tilde w_t = V^{\top}\tilde v_t$. Multiplying (\ref{coupled_dynamics}) from the left by $V^{\top}\otimes I_n$, we obtain
\begin{align}
\tilde{z}_{t+1} = \left( aI_N - g\left( V^{\top}({\cal L}+{\cal L}_R)V \right)\right ) \tilde{z}_t - \tilde{\psi} \left (\tilde z_t\right ) + \tilde w_t, \label{coupled_dyn_trans}
\end{align}
where $\tilde{\psi}(\tilde z_t) = V^{\top}\tilde{\phi}\left(\tilde x_t\right)$. We can now write $\tilde z_t = \left[\begin{array}{cc}\bar{x}_t & \hat{z}_t^{\top} \end{array}\right]^{\top}$, $\tilde{\psi}(\tilde z_t) := \left[\begin{array}{cc} \bar{\phi}_t & \hat{\psi}_t^{\top} \end{array}\right]^{\top}$,
and
$\tilde{w}_t := \left[\begin{array}{cc}
\bar{v}_t^{\top} & \hat{w}_t^{\top}
\end{array}\right]^{\top}$, where $\bar{x}_t := \frac{\mathds{1}^{\top}}{\sqrt{N}}\tilde{x}_t = \frac{1}{\sqrt{N}}\sum_{k = 1}^N x^k_t$, $\hat{z}_t := U^{\top}\tilde{x}_t$, $\bar{\phi}_t := \frac{\mathds{1}^{\top}}{\sqrt{N}}\tilde{\phi}\left(\tilde z_t\right) = \frac{1}{\sqrt{N}}\sum_{k = 1}^N \phi(x^k_t)$, $\hat{\psi}_t := U^{\top}\tilde{\phi}\left(\tilde x_t\right)$, $\bar{v}_t := \frac{\mathds{1}^{\top}}{\sqrt{N}}\tilde{v}_t = \frac{1}{\sqrt{N}}\sum_{k = 1}^N v^k_t$, and $\hat{w}_t := U^{\top}\tilde{v}_t$. Furthermore, we have $E_{\tilde v}[\bar{v}_t^2] = \sqrt{N}\omega^2,\; E_{\tilde v}[\hat{w}_t\hat{w}_t^{\top}] = U^{\top}E_{\tilde v}[\tilde{v}_t\tilde{v}_t^{\top}]U = \omega^2I_{N-1}$ and ${\cal L}_R = \sum_{e_{ij}\in \Ed_U}\xi_{ij}\ell_{ij}\ell_{ij}^{\top}$, where $\ell_{ij} \in \mathbb{R}$ is $1$ and $-1$ in the $i^{th}$ and $j^{th}$ entries, respectively, and zero elsewhere. Thus, $\ell_{ij}^{\top}\ell_{ij} = 2$ for all $e_{ij} \in \Ed$. Hence, if $\hat{\ell}_{ij} = U^{\top}\ell_{ij}$, we have $\hat{\ell}_{ij}^{\top}\hat{\ell}_{ij} = 2$ for all edges $e_{ij} \in \Ed_U$ and $U^{\top}{\cal L}_RU = \sum_{e_{ij}\in \Ed_U}\xi_{ij}\hat{\ell}_{ij}\hat{\ell}_{ij}^{\top}$. From \eqref{coupled_dyn_trans}, we obtain $\bar{x}_{t+1} =  a \bar{x}_t - \bar{\phi}_t + \bar{v}_t$ and
\begin{align}\label{sync_dynamics}
\hat{z}_{t+1} &= \left( aI_{N-1} -  g \hat{\Lambda} + g\sum_{e_{ij}\in \Ed_U}\xi_{ij}\hat{\ell}_{ij}\hat{\ell}_{ij}^{\top} \right) \hat{z}_t - \hat{\psi}_t + \hat{w}_t := A(\Xi) \hat{z}_t - \hat{\psi}_t + \hat{w}_t,
\end{align}
where $\hat{\Lambda} = \text{diag}\{\lambda_2,\cdots,\lambda_N\}$ and $\Xi = \{\xi_{ij}\;|\;e_{ij}\in \Ed_U\}$. For the synchronization of system (\ref{coupled_dynamics}), we only need to demonstrate the mean square stability about the origin of the $\hat{z}$ dynamics as given in (\ref{sync_dynamics}).

The objective is to synchronize, in a mean square sense, $N$ first-order systems over a network with a nominal graph Laplacian ${\cal L}$ with eigenvalues $0=\lambda_1 < \lambda_2 \leq \cdots \leq \lambda_N$ and maximum link CoD $\bar{\gamma}$. We present the main result of this paper.

\noindent {\bf Mean Square Synchronization Result}: The network system in Eq. \eqref{coupled_dynamics} is MSS if there exists a positive constant $p < \delta$ that satisfies
\begin{align}\label{sync_ric0}
\left( p - \frac{1}{\delta} \right) \left( \frac{1}{p} - \frac{1}{\delta} \right) > \alpha_0^2,
\end{align}
where $\alpha_0^2 = (a_0 - \lambda_{sup} g)^2 + 2\bar{\gamma}\tau\lambda_{sup}g^2$, $a_0 = a - \frac{1}{\delta}$ and $\lambda_{sup} = \underset{\lambda \in \{\lambda_2,\lambda_N\}}{\text{argmax}}\Big|\lambda+\bar{\gamma}\tau-\frac{a_0}{g}\Big|$.
Furthermore, $\tau := \frac{\lambda_{N_U}}{\lambda_{N_U}+\lambda_{2_D}}$, where $\lambda_{N_U}$ is the maximum eigenvalue of ${\cal L}_U$ and $\lambda_{2_D}$ is the second-smallest eigenvalue of ${\cal L}_D$.

The derivation of this result will be discussed in the Methods section. The above synchronization result relies on a Lyapunov function-based stability theorem. The positive constant $p$ in Eq. (\ref{sync_ric0}) is used in the construction of the Lyapunov function given by $V(x_t)=px_t^2$. Furthermore, in the Methods section, we prove that the Mean Square Synchronization Result obtained in \eqref{sync_ric0} is equivalent to
\begin{align}\label{sync_ric00}
\left(1 - \frac{1}{\delta}\right)^2 > \alpha_0^2 = \left(a_0 - \lambda_{sup}g\right)^2 + 2\bar{\gamma}\tau\lambda_{sup}g^2.
\end{align}
The main result can be interpreted in multiple ways. One particular interpretation useful in the subsequent definition of the synchronization margin is adapted from robust control theory. The robust control theory results allow one to analyse the stability of the feedback system with uncertainty in the feedback loop. The basic concept is that if the product of the system gain and the gain of the uncertainty (also called the loop gain) are less than one, then the feedback system is stable \cite{Astrom_book}. Note that system and uncertainty gains are measured by appropriate norms. The farther the system gain is from unity, the more uncertainty the feedback loop can tolerate and hence the more robust the system is to uncertainty. This result from robust control theory is extended to the case of stochastic uncertainty and nonlinear system dynamics \cite{amit_synchronization, amit_erasure_observation_journal,amit_ltv_journal,Vaidya_erasure_SCL,Vaidya_erasure_stablization, scl04, Haddad_book}. It can be shown that the synchronization problem for network systems with stochastic uncertainty can be written in this robust control form, where the loop gain directly translates to the synchronization margin. We refer the reader to supplementary material for more details and a mathematically rigorous discussion on the robust control-based interpretation behind the following mean square synchronization margin definition.

\noindent {\bf Mean Square Synchronization Margin:} The equivalent Mean Square Synchronization Result is used to define the Mean Square Synchronization Margin as follows:
\begin{align}\label{sync_margin}
\rho_{SM} := 1 - \frac{\sigma^2g^2}{\left(1 - \frac{1}{\delta}\right)^2 - \hat{a}^2}
\end{align}
where $\hat{a} = a - \frac{1}{\delta} - \mu g$, $\hat{a}^2 < \left(1 - \frac{1}{\delta}\right)^2$, $\mu = \lambda_{sup}$, $\sigma^2 = 2 \bar{\gamma}\tau\lambda_{sup}$, and $\lambda_{sup} := \underset{\lambda \in \{\lambda_2,\lambda_N\}}{argmax}\Big|\lambda+\bar{\gamma}\tau-\frac{a_0}{g}\Big|$. Furthermore, $\tau := \frac{\lambda_{N_U}}{\lambda_{N_U}+\lambda_{2_D}}$, where $\lambda_{N_U}$ is the maximum eigenvalue of ${\cal L}_U$ and $\lambda_{2_D}$ is the second-smallest eigenvalue of ${\cal L}_D$.

$\rho_{SM}$ measures the degree of robustness to stochastic perturbation. In particular, the larger the value of $\rho_{SM}$ (i.e., the smaller the value of $\left(\frac{1}{\left(1 - \frac{1}{\delta}\right)^2 - \hat{a}^2}\right)$), the larger the variance of stochastic uncertainty that can be tolerated in the network interactions before the network loses synchronization. When considering practical computation, it is important to emphasize $\rho_{SM}$, as computed by Eq. \eqref{sync_margin}, is obtained from a sufficiency condition and hence is a guaranteed synchronization margin, i.e., the true synchronization margin will be larger than or equal to $\rho_{SM}$. The synchronization condition for MSS of an $N$-node network system \eqref{coupled_dynamics} as formulated in Eq. \eqref{sync_margin} is provided in terms of a scalar quantity instead of an $N$-dimensional matrix inequality. The condition is independent of network size, which makes it computationally attractive for large-scale networks. We now discuss the effects of various network parameters on synchronization.

\noindent {\bf Role of $\tau$ and $\bar{\gamma}$:} The parameter $0 < \tau \leq 1$ in $\rho_{SM}$ captures the effect of the uncertainty location in the graph topology. If the number of uncertain links ($|\Ed_U|$) is large, the deterministic graph will become disconnected ($\lambda_{2_D} = 0$), and thus $\tau$ will equal $1$. In contrast, if a single link is uncertain ($\Ed_U = \{e_{kl}\}$), then $\tau = \frac{2\mu_{kl}}{2\mu_{kl}+\lambda_{2_D}}$. This indicates that the synchronization degradation is proportional to the link weight. Because $\lambda_{2_D} \leq \lambda_2$, a lower algebraic connectivity of the deterministic graph further degrades $\rho_{SM}$. Thus, we can rank-order individual links within a graph with respect to their degradation of $\rho_{SM}$, where a smaller $\tau$ produces an increased $\rho_{SM}$.
 For example, it can be proved that the average value of $\tau$ for a nearest neighbour network is larger than that for a random network \cite{Nature_watts}. Thus, if a randomly chosen link is made stochastic in a nearest neighbour network and in a random network, the margin of synchronization decreases by a larger amount in the nearest neighbour network as compared than in the random network. We provide simulation results to support this claim in the supplementary information section. The significance of $\bar{\gamma}$ is straightforward, as it captures the maximum tolerable variance of the system, normalized with respect to the mean weight of the link. If $\bar{\gamma} > 1$, then the uncertainty occurring within the system is clustered, which leads to large intervals of high deviation. Similarly, if $\bar{\gamma} < 1$, then the uncertainties are bundled closer to the mean value. Decreasing $\bar{\gamma}$ for the network increases $\rho_{SM}$.

\noindent {\bf Role of Laplacian Eigenvalues:} The second smallest eigenvalue of the nominal graph Laplacian $\lambda_2>0$ indicates the algebraic connectivity of the graph. Because $\alpha_0$ in \eqref{sync_margin} is a quadratic in $\lambda$, there exist critical values of $\lambda_2$ (or $\lambda_N$) for a given set of system parameters and CoD below which (or above which) synchronization is not guaranteed. Hence, the critical $\lambda_2$ indicates that there is a required minimum degree of connectivity within the network for synchronization to occur. Furthermore, increasing the connectivity at appropriate nodes may increase $\lambda_2$, leading to higher $\rho_{SM}$. To understand the significance of $\lambda_N$, we look at the complement of the graph on the same set of nodes. We know from \cite{Merris1994} (Lemma provided in Supplementary Information for reference) that the sum of the largest Laplacian eigenvalue of a graph and the second smallest Laplacian eigenvalue of the complementary graph is a constant. Thus, if $\lambda_N$ is large, then the complementary graph has low algebraic connectivity. Hence, a high $\lambda_N$ indicates the presence of many densely connected nodes. Therefore, we conclude that a robust synchronization is guaranteed for graphs with close-to-average node connectivity to graphs with isolated but highly connected hub nodes. Thus, decreasing $\lambda_N$ by reducing the connectivity of specific nodes (i.e., dense hub nodes) will help increase $\rho_{SM}$.

\noindent {\bf Impact of Internal Dynamics:} The internal dynamics are captured by parameters $a$ and $\delta$, which respectively represent the rate of linear instability and the bound on the rate of change of the nonlinearity. As $a$ increases, the linear dynamics become more unstable. When all other parameters are held constant, an increase in $a$ results in a decrease in $\left(1-\frac{1}{\delta}\right) - \left(a_0 - \lambda_{sup}g\right)^2$. Because $\rho_{SM} \propto -\frac{1}{\left(1 - \frac{1}{\delta}\right)^2 - \left(a_0 - \lambda_{sup}g\right)^2}$, an increase in $a$ will produce a decrease in $\rho_{SM}$. Thus, as the instability of the internal dynamics increases, the network becomes less robust to uncertainty. When the fluctuations in link weights are zero (i.e., CoD $\bar{\gamma} = 0$), the critical value of $\lambda_2$ below which synchronization is not guaranteed is $\lambda_2^* = \frac{a-1}{g}$. Furthermore, synchronization is not guaranteed for $\lambda_N$ above the critical value $\lambda_N^* = \frac{a+1}{g} - \frac{2}{g\delta} = \lambda_2^* + \frac{2}{g}\left(1 - \frac{1}{\delta}\right)$. Thus, we see $\lambda_N^* - \lambda_2^* = \frac{2}{g}\left(1 - \frac{1}{\delta}\right)$ and $\frac{\lambda_N^*}{\lambda_2^*} = 1 + \frac{2}{a-1}\left(1 - \frac{1}{\delta}\right)$. While $\lambda_N^* - \lambda_2^*$ 
is independent of the internal dynamics parameter $a$, $\lambda_2^*$ increases with an increase in $a$. In fact, for $a=1+\epsilon$, where $\epsilon>0$ is arbitrarily small, we have $\lambda_2^*=\frac{\epsilon}{g}$. Hence, as the internal dynamics become more unstable, we require a higher degree of connectivity between the network agents to achieve synchronization. Because the nonlinearity $\phi$ is sector-bounded by $\frac{2}{\delta}$, the impact of the nonlinearity on synchronization can be analysed using $\delta$. When all of the other network parameters are held constant, $\lambda_2^*$ is independent of $\delta$ and $\lambda_N^*$ increases with increasing $\delta$. Increasing the value of $\delta$ leads to an increase in $\left(1-\frac{1}{\delta}\right)^2 - \left(a_0 - \lambda_{sup}g\right)$, which increases $\rho_{SM}$. Hence, as the nonlinearity of the system is reduced, the system becomes more robust to uncertainties.

\noindent {\bf Impact of Coupling Gain:} The impact of the coupling gain is more complicated than the impact of the internal dynamics. A very small coupling gain is not enough to guarantee $\left(1-\frac{1}{\delta}\right)^2 > (a_0-\lambda_{sup}g)^2 + 2\bar{\gamma}\tau g^2$, which is required to ensure $\rho_{SM} > 0$. On the other hand, a very large coupling gain also does not guarantee $\left(1-\frac{1}{\delta}\right)^2 > (a-\lambda_{sup}g)^2 + 2\bar{\gamma}\tau g^2$. Thus, we can conclude the coupling gain affects the synchronization margin in a nonlinear fashion. Hence, to obtain the largest possible $\rho_{SM}$, the network must operate at an optimal gain.

We now demonstrate how the main results of this paper can be used to determine the optimal value of the coupling gain $g^*$ that maximizes the margin of synchronization for a given network topology (i.e., specific values of $\lambda_2$ and $\lambda_N$) and uncertainty (i.e., CoD value $\bar \gamma$). We assume that, for given values of $\lambda_2, \lambda_N$, and $\bar \gamma$, there exists a value of $g$ for which synchronization is possible.

\noindent {\bf Optimal Gain:} For the network system in Eq. \eqref{coupled_dynamics} with $\rho_{SM}$ given by Eq. \eqref{sync_margin}, the optimal gain $g^*$ that produces the maximum $\rho_{SM}$ is
\begin{align}\label{opt_gain}
g^* = \frac{2(a-\frac{1}{\delta})}{\max\{\lambda_N,\lambda_2+2\bar{\gamma}\tau \}+\lambda_2 + 2\bar{\gamma}\tau}.
\end{align}

The derivation of this result will be discussed in the Methods Section. The results of the Mean Square Synchronization Margin $\rho_{SM}$ and the Optimal Gain $g^*$ will be used in the following subsections to study the effect of neighbours and network connectivity on both nearest neighbour networks and random networks such as Erdos-Renyi and Small-World networks.

\subsection{Interplay of Internal Dynamics, Network Topology, and Uncertainty Characteristics}\label{interplay}

We now study the interplay of the internal dynamics ($a$), nonlinearity bound ($\delta$), network topology ($\lambda$), and the uncertainty characteristics ($\bar{\gamma}$) through simulations over a $1000$-node network using a set of parameter values. To nullify the bias of uncertain link locations, we choose to work with a large number of uncertain links to obtain $\tau \approx 1$.

\begin{figure*}[ht]
\begin{center}
\includegraphics[width = 1\textwidth]{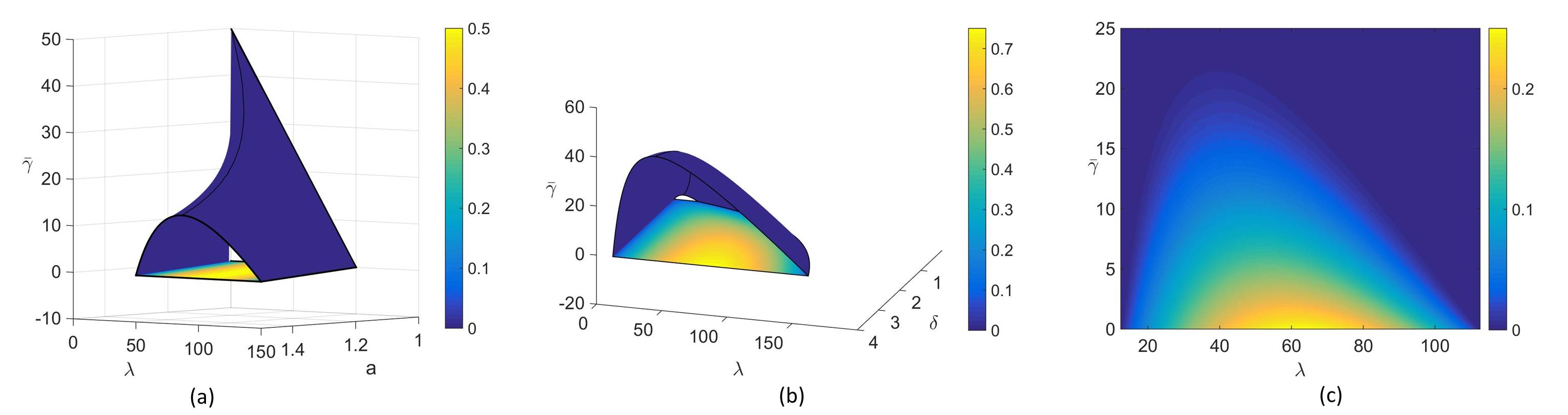}
\caption{(a) $\rho_{SM}$ in $a-\lambda-\bar \gamma$ parameter space for $g=0.01$ and $\delta=2$, (b) $\rho_{SM}$ in $\delta-\lambda-\bar \gamma$ parameter space for $a=1.125$ and $g=0.01$, (c) $\lambda-\bar \gamma$ parameter space indicating $\rho_{SM}$ for $a=1.125$, $g=0.01$, and $\delta=2$.}
\label{fig1s2}
\vspace{-0.5pc}
\end{center}
\end{figure*}

In Fig. \ref{fig1s2}(a), we study the interplay of network topology, uncertainty, and the internal dynamics in the three-dimensional parameter space  of $a-\lambda-\bar \gamma$. In Fig. \ref{fig1s2}(a), the region inside (or outside) the tunnel corresponds to the combination of parameter values where synchronization is possible (or not possible).  Another important observation we make from Fig. \ref{fig1s2}(a) is that the area inside the tunnel increases with a decrease in either the internal instability or $a$.
In Fig. \ref{fig1s2}(b), we plot the effects of changing the nonlinearity bound $\delta$ on the synchronization margin in the $\delta-\lambda-\bar{\gamma}$ space. As $\delta$ is increased, the region of synchronization increases. Thus, a minimally nonlinear system is able to achieve synchronization even with high levels of communication. On the other hand, as the nonlinearity in a system becomes significant, the interaction between the nonlinearity and the fluctuations in the link weights could have adverse effects in a highly connected network. Intuitively, because a high communication amplifies the uncertainty between the agents, one might view this as the uncertainty in the fluctuations being wrapped around and amplified by the nonlinearity, which causes this high-communication desynchronization . In Fig. \ref{fig1s2}(c), we plot a slice of the synchronization regions from both Figs. \ref{fig1s2}(a) and \ref{fig1s2}(b) for $a=1.125$, $\delta = 2$, and $g = 0.01$, that highlights the synchronization margin.

\subsection{ Optimal Neighbours in Nearest Neighbour Networks}\label{neighbour}
\begin{figure*}[htb!]
\begin{center}
\includegraphics[width = 1\textwidth]{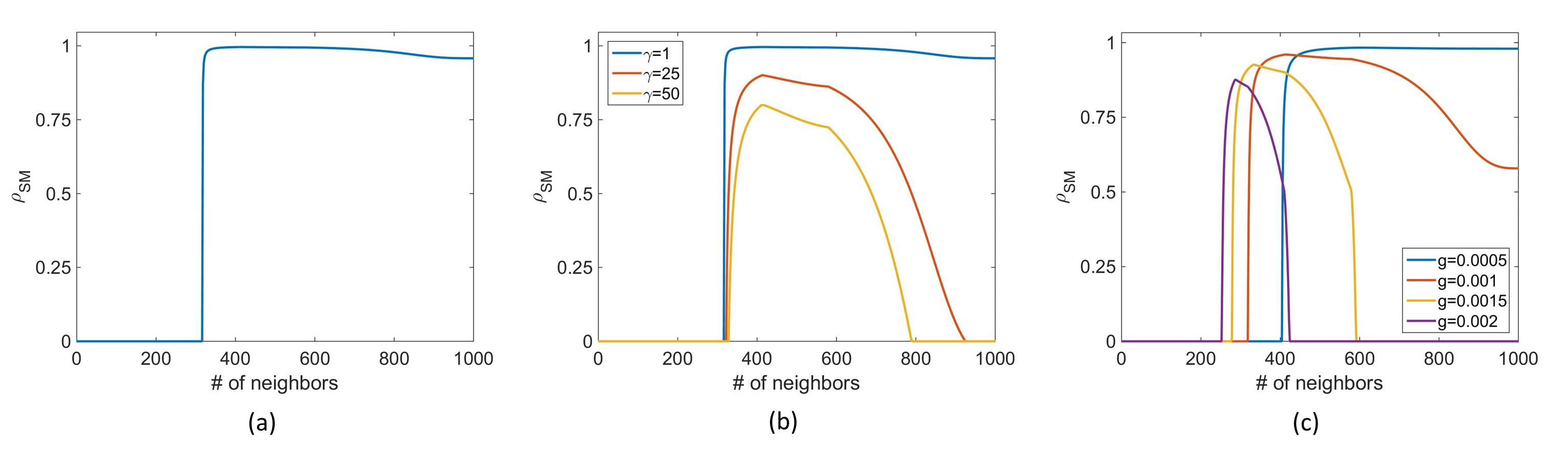}
\caption{(a) Synchronization margin for $a = 1.05$, $\delta = 2$, $g=0.001$, and $\bar{\gamma} = 1$ as the number of neighbours are varied in a nearest neighbour graph, (b) Synchronization margin for $a = 1.05$, $\delta = 2$, and $g = 0.001$ for different $\bar{\gamma}$ as the number of neighbours are varied in a nearest neighbour graph, (c) Synchronization margin for $a = 1.05$, $\delta = 2$, and $\bar{\gamma} = 10$ for different coupling gains as the number of neighbours are varied in a nearest neighbour graph.}
\label{NNstability}
\end{center}
\end{figure*}

The analytical formula for the synchronization margin in Eq. \eqref{sync_margin} provides us with a powerful tool to understand the effect of various network parameters on the synchronization margin. In this section, we investigate the effects of the number of neighbours on the synchronization margin. We consider a nearest neighbour network with $N=1000$ nodes and increase the number of neighbours to study their impact on the synchronization margin. The other network parameters are set to $a=1.05, \delta=2$, $g=\frac{1}{N}$, and $\bar{\gamma}=25$. We choose a large number of uncertain links (70\%) so that $\tau \approx 1$ to remove the bias of uncertain link locations. We show the plot for the synchronization margin versus the number of neighbours in Fig. \ref{NNstability}(a). From this plot, we see that there exists an optimal number of neighbours an agent requires in order to maximize the synchronization margin. Additionally, there is a minimum number of neighbours required by any given agent. Below this number, the network will not synchronize. However, an uncertain environment with too many neighbours is also detrimental to synchronization. This result highlights the fact that, while ``good'' information is propagated through neighbours via network interconnection, in an uncertain environment, these same neighbours can propagate ``bad'' information that is detrimental to reaching an agreement. In Fig. \ref{NNstability}(b), we show the plot for the change in the synchronization margin versus a change in the number of neighbours for different values of CoD. We notice that, for larger values of CoD, the drop in the margin as the network connectivity increases is more dramatic.

In light of the previous discussion, we can also interpret the coupling gain $g$ as the amount of trust a given agent has in the information provided by its neighbours. In particular, if the coupling gain is large, then the agent has more trust in its neighbours. In Fig. \ref{NNstability}(c), we show the effects of increasing the coupling gain on the synchronization margin. We observe that if an agent has more trust in its neighbours, then fewer neighbours are required to achieve synchronization. However, in an uncertain environment, an agent with more trust in its neighbours must avoid having more neighbours, as it is detrimental to synchronization. On the other hand, if an agent has less trust in its neighbours, more connections must be formed to gather as much information as possible, even if that information is corrupted. Thus, forging connections is good for a group with the goal of synchronization, but there exists a critical number of neighbours above which the benefits from forging new connections diminish.

\subsection{Optimal gain for complex networks}

Based on the optimal gain formulation, we can now compare the performance of some well-known random networks and the optimal gain required to synchronize these networks. We use the following parameters in these simulations: the system instability $a=1.05$, the nonlinearity bound $\delta = 4$, and the uncertainty statistics represented by CoD is $\bar{\gamma}=1$. Furthermore, we choose $\tau \approx 1$. The properties of these random networks are studied for four different network sizes: $N \in \{80,100,120,140\}$, where $N$ is the number of nodes.

\begin{figure*}[ht]
\begin{center}
\includegraphics[width = 0.8\textwidth]{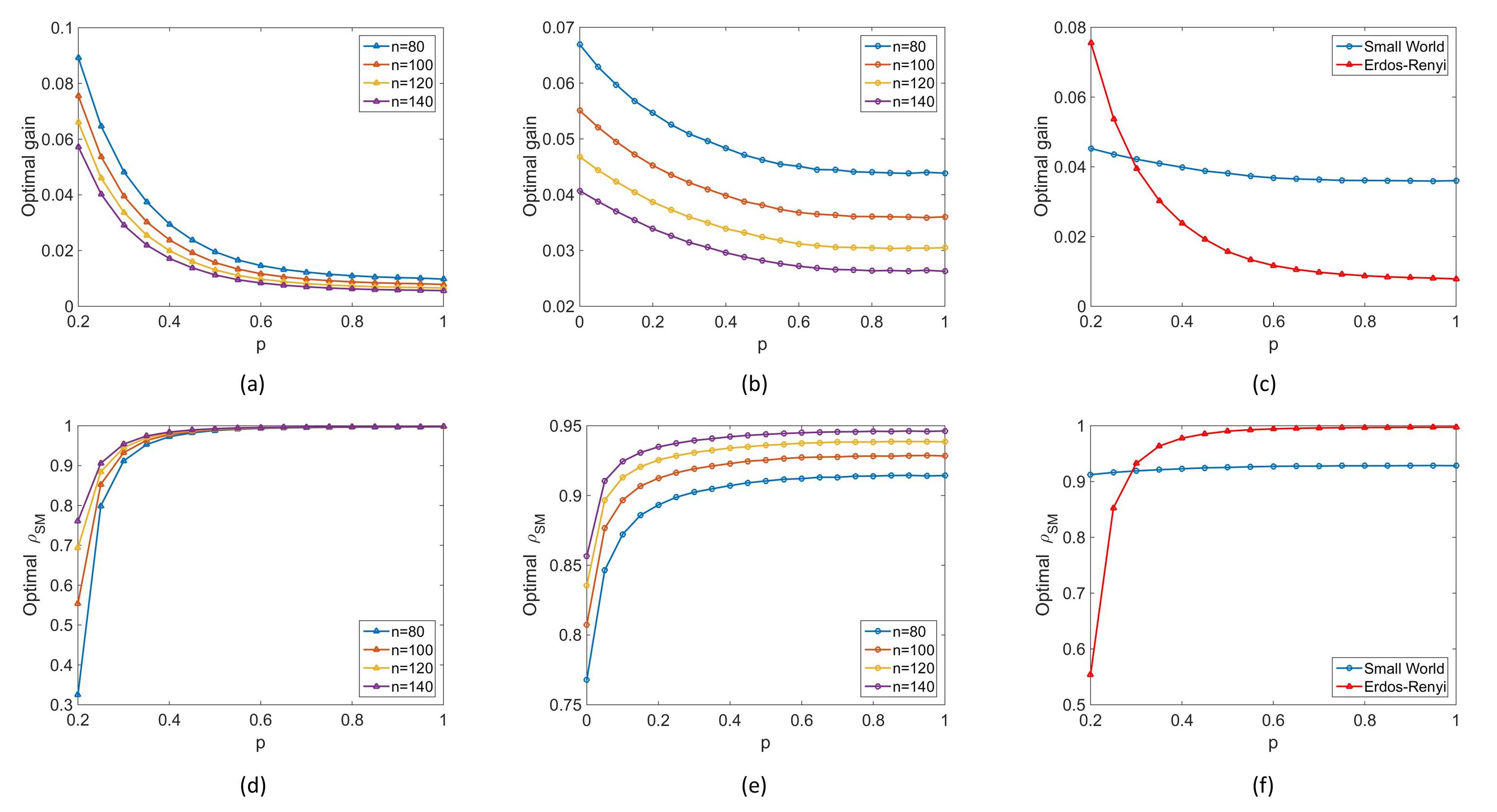}
\caption{ Optimal gain computation for (a) an Erdos-Renyi network with probability of connecting two nodes $p$, for varying network sizes and (b) a Small World network with probability of rewiring an edge $p$, for varying network sizes; (c) comparison of optimal gain for Erdos-Renyi and Small World networks as a function of probability for network size $n=100$. Optimal synchronization margin computation for (d) Erdos-Renyi network with probability of connecting two nodes $p$, for varying network sizes and (e) a Small World network with probability of rewiring an edge $p$, for varying network sizes; (f) comparison of optimal synchronization margin for Erdos-Renyi and Small World networks as a function of probability for network size $n=100$.}
\label{fig1s5}
\end{center}
\end{figure*}

In Fig. \ref{fig1s5}(a), we plot the optimal gain for the Erdos-Renyi (ER) networks as a function of the edge connection probability. It is well known that for an Erdos-Renyi network of size $N$ to be connected, the probability of connection must be $p \geq \frac{\log{N}}{N}$. Hence, we plot these networks for probabilities ranging from $p=0.2$ to $p=1$. At $p=1$, we obtain an all-to-all connection network, as each edge is connected with unit probability. In Fig. \ref{fig1s5}(d), we plot the corresponding optimal synchronization margin for the ER network. In Figs. \ref{fig1s5}(b) and \ref{fig1s5}(e), we plot the optimal gain and optimal synchronization margin, respectively, for a SW network with varying probability $p$ \cite{Nature_watts}. To better observe the contrast in behaviour of both the ER and SW random networks, we plot in Fig. \ref{fig1s5}(c) the optimal gains for an ER network and an SW network with $N=100$ nodes.

We notice that, while a larger gain is required to synchronize the ER network than that for the SW network for smaller values of $p$, the optimal gain for the ER network is smaller than that of the SW network for larger values of $p$. In Figs. \ref{fig1s5}(d), \ref{fig1s5}(e), and \ref{fig1s5}(f), we plot the optimal synchronization margins for the two networks. We notice an increase in the synchronization margin for the ER network around $p=0.5$. From these plots (specifically Figs. \ref{fig1s5}(c) and \ref{fig1s5}(f)), we conclude that for the given set of parameters, the ER (or SW) network has better synchronization properties (i.e., a smaller value of the optimal gain and a larger margin of synchronization) for larger (or smaller) values of $p$. The transition between the two cases occurs for some probability between $p=0.2$ and $p=0.4$.

\section{Discussion}\label{section_conc}
We study the problem of synchronization in complex network systems in the presence of stochastic interaction uncertainty between the network nodes. We exploited the identical nature of the internal node dynamics to provide a sufficient condition for network synchronization. The unique feature of this sufficient condition is its independence from the network size. This makes the sufficient condition computationally attractive for large-scale network systems. Furthermore, this sufficient condition provides useful insight into the interplay between the internal dynamics of the network nodes, the network interconnection topology, the location of uncertainty, and the statistics of the uncertainty and into their effects on the network synchronization. The sufficient condition provided in the main result allows us to characterize the degree of robustness of a synchronized state to stochastic uncertainty through the definition of a mean square synchronization margin. Using the synchronization margin, a formulation for an optimal synchronization gain is derived to assist in designing gains for complex networks based purely on the system dynamics, nominal network Laplacian eigenvalues, and uncertainty statistics. This optimal gain result is used to compare various complex network topologies for given internal nodal dynamics.

When considered from a practical point of view, the synchronization margin is useful in determining the synchronizability of large-scale networks with stochastic uncertainty in the coupling. The independence of the result with respect to the network size can be used to obtain a bound on the tolerable uncertainty with minimal computational effort. In networked systems with communication uncertainty, these results can be used to provide a worst-case signal-to-noise ratio that is tolerable in communication or to design network connectivity in order to optimize the network's tolerance to uncertainty. These results have potential applications in determining the optimal neighbours and coupling gain in consensus dynamics, swarm dynamics, and other situations where systems seek synchronization.


\section{Methods}

\noindent {\bf Mean Square Synchronization Condition:} The system described by Eq. \eqref{coupled_dynamics} is MSS as given by Definition 1, if there exist $L > 0$, $K > 0$ and $0 < \beta < 1$, such that
\begin{align}\label{z_MSS}
E_{\Xi_0^{t-1},v_0^{t-1}} \left[ \parallel \hat z_t \parallel^2 \right] &\le   K  {\beta}^t \parallel \hat z_{0}  \parallel^2 + L\omega^2.
\end{align}
We refer to this as mean square stability of $\hat{z}_t$. From Eq. \eqref{sync_dynamics}, we obtain, $\parallel \hat z_t \parallel^2 = \tilde x_t^{\top} \left(UU^{\top} \otimes I_n\right)\tilde x_t = \frac{1}{2N}\sum_{i=1}^N\sum_{j\neq i,j=1}^N \big\|x_t^i - x_t^j\big\|^{2}$, since $UU^{\top} = I_N - \frac{1}{N}\bf{1}\bf{1}^{\top}$.
Now, suppose there exist $L > 0$, $K > 0$, and $0 < \beta < 1$, such that \eqref{z_MSS} holds true. We can rewrite \eqref{z_MSS} as
\begin{align}\label{sum_error}
&E_{\Xi_0^{t-1},v_0^{t-1}} \left[\sum_{k=1}^N\sum_{j\neq k, j=1}^N \big\| x^k_t - x^j_t \big\|^2 \right] \le  K  {\beta}^t  \sum_{k=1}^N\sum_{j\neq k,j=1}^N \parallel x^k_0 - x^j_0  \parallel^2 + 2NL\omega^2.
\end{align}
Thus, from \eqref{sum_error} we obtain systems $S_k$ and $S_l$, that satisfy \eqref{MSS_def} for mean square synchronization, where $\bar K(\tilde{e}_0) := K \left ( 1 + \frac{\sum_{i=1, i\neq k}^N\sum_{j=1, j\neq i}^N \parallel x^i_0 - x^j_0 \parallel^2}{\parallel x^k_0 - x^l_0 \parallel^2}  \right )$ and $\bar{L}=2NL$.

In the Mean Square Synchronization Condition, we proved the mean square stability of \eqref{sync_dynamics} guarantees the MSS of \eqref{coupled_dynamics}. We will now utilize this result to provide a sufficiency condition for MSS of \eqref{sync_dynamics}.\\

\noindent {\bf Mean Square Stability of the Reduced System:} The system given by \eqref{sync_dynamics} is mean square stable, if there exists a Lyapunov function $V(\hat{z}_t) = \hat{z}_t^{\top}P\hat{z}_t$ for a symmetric matrix $P > 0$, such that for some symmetric matrix $R_P >0$ and $\rho > 0$ we have,
\begin{align}\label{1step_lyap}
E_{\Xi_t,v_t} \left  [ V(\hat{z}_{t+1})) - V(\hat{z}_t) \right ] < - \hat{z}_t^{\top}R_P\hat{z}_t + \rho \omega^2.
\end{align}

Consider $V(\hat{z}_t) = \hat{z}_t^{\top}P\hat{z}_t$ for a symmetrix matrix $P>0$, we know there exist $0< c_1 < c_2$, such that $c_1\|\hat{z}_t\|^2 \leq V_t \leq c_2 \|\hat{z}_t\|^2$.
Let $V(\hat{z}_t)$ satisfy \eqref{1step_lyap}. Substituting $c_3 = \lambda_{max}(R_P)$ as the spectral radius of $R_P$ in \eqref{1step_lyap} and using $c_2$ sufficiently large to define $\beta := 1 - \frac{c_3}{c_2} > 0$, we obtain, $E_{\Xi_t,v_t} \left  [ V(\hat{z}_{t+1})) \right ] < \beta V(\hat{z}_t) + \rho \omega^2$. Taking expectation over $(\Xi_0^t,v_0^t)$ recursively, we obtain, $c_1 E_{\Xi_0^t,v_0^t} \left[ \|\hat{z}_{t+1}\|^2 \right] < c_2 \beta^{t+1} \parallel \hat{z}_0\parallel^2 + \frac{1}{1-\beta}\rho \omega^2$. This guarantees the mean square stability of $\hat{z}_t$, for $K = \frac{c_2}{c_1}$ and $L = \frac{\rho}{(1-\beta)c_1}$.

We now utilize the Mean Square Stability of the Reduced System to define the Mean Square Synchronization Margin as given in \eqref{sync_margin}. Towards this aim, we first construct an appropriate Lyapunov function, $V(\hat{z}_t) = \hat{z}_t^{\top}P\hat{z}_t$, that guarantees mean square stability. From \eqref{sync_dynamics}, defining $\Delta V :=E_{\Xi_t,v_t}\left[ V(\hat{z}_{t+1}) - V(\hat{z}_t) \right]$, we obtain,
\begin{align}
\Delta V &= E_{\Xi_t}\left[\hat{z}_t^{\top}\left(A(\Xi_t)^{\top}PA(\Xi_t) - P\right)\hat{z}_t -\hat{z}_t^{\top}A(\Xi_t)^{\top}P\hat{\psi}_t - \hat{\psi}_t^{\top}PA(\Xi_t)\hat{z}_t + \hat{\psi}_t^{\top}P\hat{\psi}_t\right]
 + E_{v_t}[\hat{w}_t^{\top}P\hat{w}_t].
\end{align}
Now, suppose for some $R_P > 0$, $P$ satisfies,
\begin{align}\label{unc_lure_p1}
P =& E_{\Xi} \left [ A(\Xi)^{\top}PA(\Xi)\right]  + R_P + E_{\Xi} \left [\left(A(\Xi)^{\top}P-I_{N-1}\right)(\delta I_{N-1}-P)^{-1}\left(P A(\Xi)-I_{N-1}\right)\right ].
\end{align}
Using \eqref{unc_lure_p1} and algebraic manipulations as given in \cite{haddad_bernstein_DT}, we can rewrite $\Delta V = - \hat{z}^{\top}_t R_P \hat{z}_t - E_{\Xi_t}\left[\eta^{\top}_t \eta_t \right] - 2 \hat\psi^{\top}_t\left(\hat{z}_t - \frac{\delta}{2}\hat\psi_t\right) + \textrm{trace}(P E_{v_t}[\hat{w}_t\hat{w}_t^{\top}])$, where $\eta_t(\Xi(t))$ is given by $\eta_t(\Xi(t)) = W^{-\frac{1}{2}}\left(PA(\Xi) - I_{N-1}\right)\hat{z}_t - W^{\frac{1}{2}}\hat\psi_t$ and $W := (\delta I_{N-1} - P)$. Since, $\phi(\cdot)$ is monotonic and globally Lipschitz with constant $\frac{2}{\delta}$, we know $\left(\phi\left(x^k_t\right)-\phi\left(x^l_t\right)\right)^{\top}\left(\frac{2}{\delta}\left(x^k_t-x^l_t\right) - \left(\phi\left(x^k_t\right)-\phi\left(x^l_t\right)\right)\right) > 0$. This gives $\hat{\psi}_t^{\top}\left(\hat{z}_t - \frac{\delta}{2}\hat{\psi}_t\right)>0$. Using this and writing $\rho = \textrm{trace}(P)$, we obtain Eq. \eqref{1step_lyap}. Hence, \eqref{unc_lure_p1} is sufficient for MSS of \eqref{component_dynamics} from condition for Mean Square Stability of the Reduced System. Furthermore, the Eq. in (\ref{unc_lure_p1}) can be rewritten using ~\cite{lancaster} (Proposition 12.1,1) as
\begin{align} \label{unc_lure_p}
P =& E_{\Xi} \left [ A_0(\Xi)^{\top}PA_0(\Xi)\right]  + R_P + \frac{1}{\delta}I_{N-1}+E_{\Xi} \left [A_0(\Xi)^{\top}P(\delta I_{N-1}-P)^{-1}P A_0(\Xi)\right ],
\end{align}
where $A_0(\Xi) = a_0I_{N-1} -g\hat{\Lambda} - gU^{\top}{\cal L}_R U$ and $a_0 = a - \frac{1}{\delta}$. We observe this condition requires us to find a symmetric Lyapunov function matrix $P$ of order $\frac{N(N-1)}{2}$. We now reduce the order of computation by using network properties. For this, consider $P = p I_{N-1}$, where $p < \delta$ is a positive scalar. This gives us $\delta I_{N-1} > P$. Using this and \eqref{sync_dynamics}, we rewrite the condition in \eqref{unc_lure_p} as follows,
\begin{align} \label{ric1}
pI_{N-1} &>  p(a_0I_{N-1} - g\hat{\Lambda})^{\top}(a_0I_{N-1}- g\hat{\Lambda})  + \frac{p^2}{\delta - p}(a_0I_{N-1}-g\hat{\Lambda})^{\top}(a_0I_{N-1}-g\hat{\Lambda}) + \frac{1}{\delta}I_{N-1}\nonumber \\
&\quad+ pg^2\underset{e_{ij \in \Ed_U}}{\sum} \sigma_{ij}^2\hat{\ell}_{ij}\hat{\ell}_{ij}^{\top}\hat{\ell}_{ij}\hat{\ell}_{ij}^{\top} +\frac{p^2g^2}{\delta - p}\underset{e_{ij \in \Ed_U}}{\sum} \sigma_{ij}^2\hat{\ell}_{ij}\hat{\ell}_{ij}^{\top}\hat{\ell}_{ij}\hat{\ell}_{ij}^{\top}.
\end{align}
We know $\hat{\ell}_{ij}^{\top}\hat{\ell}_{ij} = \ell_{ij}^{\top}U_dU_d^{\top}\ell_{ij} = \ell_{ij}^{\top}\ell_{ij} = 2$ and $\underset{e_{ij \in \Ed_U}}{\sum}\sigma_{ij}^2\hat{\ell}_{ij}\hat{\ell}_{ij}^{\top}\hat{\ell}_{ij}\hat{\ell}_{ij}^{\top}\leq 2\bar{\gamma}\underset{e_{ij \in \Ed_U}}{\sum}\mu_{ij}\hat{\ell}_{ij}\hat{\ell}_{ij}^{\top} = 2\bar{\gamma}U^{\top}{\cal L}_U U$. For $\tau = \frac{\lambda_{N_U}}{\lambda_{N_U}+\lambda_{2_D}}$, we have, ${\cal L}_U\leq \tau \left({\cal L}_D + {\cal L}_U\right) = \tau {\cal L}$. Hence, $\underset{e_{ij \in \Ed_U}}{\sum}\sigma_{ij}^2\hat{\ell}_{ij}\hat{\ell}_{ij}^{\top}\hat{\ell}_{ij}\hat{\ell}_{ij}^{\top} \leq 2\bar{\gamma}\tau U^{\top} {\cal L} U = 2\bar{\gamma}\tau \hat{\Lambda}$.
Substituting this into \eqref{ric1}, a sufficient condition for inequality \eqref{ric1} to hold is given by
\begin{align} \label{ric4}
 pI_{N-1} > \left(p+\frac{p^2}{\delta - p}\right)(a_0I_{N-1} - g\hat{\Lambda})^{\top}(a_0I_{N-1}- g\hat{\Lambda}) + \left(p+\frac{p^2}{\delta - p}\right)2\bar{\gamma}\tau g^2\hat{\Lambda} + \frac{1}{\delta}I_{N-1}.
\end{align}
Equation \eqref{ric4} is a block diagonal equation. The individual blocks provide the sufficient condition for MSS as, $p >  \left(p+\frac{p^2}{\delta - p}\right)\left((a_0 - g\lambda_j)^2 + 2\bar{\gamma}\tau g^2\lambda_j\right) + \frac{1}{\delta}$, for all eigenvalues $\lambda_j$ of $\hat{\Lambda}$. This is simplified as
\begin{align}\label{suff_main1}
\left(\frac{1}{p} - \frac{1}{\delta}\right)\left(p - \frac{1}{\delta}\right) > \alpha_0^2,
\end{align}
where $\delta > p > 0$ and $\alpha_0^2 = (a_0 - g\lambda)^2 + 2\bar{\gamma}\tau \lambda g^2$ for all $\lambda \in \{\lambda_2,\ldots,\lambda_N\}$ are eigenvalues of the nominal graph Laplacian. Now, for each of these conditions to hold true, we must satisfy condition \eqref{suff_main1} for the minimum value of $\alpha_o^2$ with respect to all possible $\lambda$. Now, $\lambda^*$ that provides minimum values for $\alpha_0^2$ is found by setting $\frac{d\alpha_0^2}{d\lambda}\Big|_{\lambda^*} = 0$, giving us $\lambda^* = \frac{a_0}{g} - \bar{\gamma}\tau$. Using $\lambda^*$, we know for \eqref{suff_main1} to be satisfied for all $\lambda \in \{\lambda_2\ldots,\lambda_N\}$, it must satisfy \eqref{suff_main1} for the farthest such $\lambda$ from $\lambda^*$. Since eigenvalues of the nominal graph Laplacian are positive and monotonic non-decreasing, all we need is to satisfy \eqref{suff_main1} for $\lambda_{sup}$, where $\lambda_{sup} = \underset{\lambda \in \{\lambda_2,\lambda_N\}}{\text{argmax}}|\lambda - \lambda^*|$.

We observe from \eqref{suff_main1}, if $p = q > 1$ is a solution of \eqref{suff_main1}, then $p = \frac{1}{q}$. We state that \eqref{suff_main1} holds, if and only if,
\begin{align}\label{suff_main2}
\left(1 - \frac{1}{\delta}\right)^2 > \alpha_0^2 = \hat{a}^2 + \sigma^2g^2,
\end{align}
where $\hat{a} = a - \frac{1}{\delta} - \mu g$, $\mu = \lambda_{sup}$, $\sigma^2 = 2\bar{\gamma}\tau \lambda_{sup}$. The ``only if'' part is obvious as, $\left(1 - \frac{1}{\delta}\right)^2 \geq \left(p - \frac{1}{\delta}\right)\left(\frac{1}{p} - \frac{1}{\delta}\right)$, from AM-GM inequality. To show the ``if'' part assume there exists $r > 0$, such that, $\left(1 - \frac{1}{\delta}\right)^2 = \alpha_0^2 + r > \alpha_0^2 + \frac{r}{2}$. Now consider some $\epsilon > 0$ such that, $\frac{r}{2} = \left(\frac{\epsilon^2}{1+\epsilon}\right)\frac{1}{\delta}$. Hence, we obtain, $\left(1+\epsilon - \frac{1}{\delta}\right)\left(\frac{1}{1 + \epsilon} - \frac{1}{\delta}\right) = 1 + \frac{1}{\delta^2} - \frac{2}{\delta} - \left(\frac{\epsilon^2}{1+\epsilon}\right)\frac{1}{\delta} = \left(1 - \frac{1}{\delta}\right)^2 - \frac{r}{2} > \alpha_0^2$. Setting $p = 1 + \epsilon > 1$, we know \eqref{suff_main1} holds true for some $p > 1$. Hence, \eqref{suff_main1} and \eqref{suff_main2} are equivalent conditions. We now use \eqref{suff_main2} to define $\rho_{SM} = 1 - \frac{\sigma^2g^2}{\left(1 - \frac{1}{\delta}\right)^2 - \hat{a}^2}$. The rationale for this and connections with existing conditions in robust control theory are discussed in the supplementary information.

We now provide the optimal coupling gain for systems with fixed internal dynamics interacting over a nominal network with a given set of uncertain links and $\bar{\gamma}$. We observe from \eqref{suff_main2}, to maximize the synchronization margin with respect to the coupling gain, $g$, we must minimize $\alpha_0^2$, with respect to $g$, and maximize $\alpha_0^2$, with respect to $\lambda$. This is a regular saddle-point optimization problem \cite{Boyd_book}. Hence, for a given $\lambda$, $\frac{\partial \alpha_0^2(\lambda,g)}{\partial g} = -2a_0\lambda + 2\left(\lambda^2 + 2\bar{\gamma}\tau\lambda\right)g = 0$. This provides us with the optimal gain as $g^*(\lambda) = \frac{a_0}{\lambda + 2\bar{\gamma}\tau}$ with $\alpha_0^2(\lambda,g^*(\lambda)) = \frac{2\bar{\gamma}\tau a_0^2}{\lambda + 2\bar{\gamma}\tau}$. The only important eigenvalues of the nominal graph Laplacian imposing limitations on synchronization, are $\lambda_2$ and $\lambda_N$. Hence we obtain $g^*(\lambda_2) = \frac{a_0}{\lambda_2 + 2\bar{\gamma}\tau}$ and $g^*(\lambda_N) = \frac{a_0}{\lambda_N + 2\bar{\gamma}\tau}$. Since $\lambda_N \geq \lambda_2$, we have $g^*(\lambda_2) \geq g^*(\lambda_N)$, and $\alpha_0^2(\lambda_2,g^*(\lambda_2)) \geq \alpha_0^2(\lambda_N,g^*(\lambda_N))$.

There also exists a value of gain, $g_e$, which provides the exact same synchronization margin for both $\lambda_2$ and $\lambda_N$. This is obtained by equating $\alpha_0^2(\lambda_2,g_e) = \alpha_0^2(\lambda_N,g_e)$, which provides, $\lambda_2^2 g_e^2 + 2\lambda_2\bar{\gamma}\tau g_e^2 - 2a_0\lambda_2 g_e = \lambda_N^2 g_e^2 + 2\lambda_N\bar{\gamma}\tau g_e^2 - 2a_0\lambda_N g_e$. This gives us, for $\lambda_N \neq \lambda_2$, and $\bar{\lambda} = \frac{\lambda_2 + \lambda_N}{2}$, $g_e = \frac{a_0}{\bar{\lambda} + \bar{\gamma}\tau}$. Furthermore, the $\alpha_0^2$ value for $g_e$, is given by, $\alpha_0^2(\lambda_2,g_e) = \alpha_0^2(\lambda_N,g_e) = a_0^2 - \frac{4\lambda_2\lambda_N a_0^2}{\left(\lambda_N + \lambda_2 + 2\bar{\gamma}\tau\right)^2}$. Since, $\lambda_N \geq \lambda_2$, we have $g_e \geq g^*(\lambda_N)$. Furthermore, $\alpha_0^2(\lambda_N,g_e) \geq \alpha_0^2(\lambda_N,g^*(\lambda_N))$, and, $\alpha_0^2(\lambda_2,g_e) \geq \alpha_0^2(\lambda_2,g^*(\lambda_2))$. We also conclude that, $g^*(\lambda_2) \geq g_e$, iff $\lambda_N \geq \lambda_2 + 2\bar{\gamma}\tau$ and $g_e \geq g^*(\lambda_2)$, iff $\lambda_2 + 2\bar{\gamma}\tau \geq \lambda_N$. We observe that, $\lambda_N \geq \lambda_2 + 2\bar{\gamma}\tau$, iff, $\alpha_0^2(\lambda_N,g^*(\lambda_2)) \geq \alpha_0^2(\lambda_N,g_e) \geq \alpha_0^2(\lambda_2,g^*(\lambda_2))$. Hence, $g_e$, being the saddle-point solution, is the optimal gain providing the largest possible $\alpha_0^2(\lambda,g)$, and the smallest $\rho_{SM}$. Similarly, $\lambda_2 + 2\bar{\gamma}\tau \geq \lambda_N$, iff, $\alpha_0^2(\lambda_2,g_e) \geq \alpha_0^2(\lambda_2,g^*(\lambda_2)) \geq \alpha_0^2(\lambda_N,g^*(\lambda_2))$. This gives $g^*(\lambda_2)$ as the optimal gain. Furthermore, at the optimal gain, we always have $\lambda_{sup} = \lambda_2$. Defining, $\chi := \max\{\lambda_N,\lambda_2+2\bar{\gamma}\tau\}$, we can write the optimal gain, $g^* = \frac{2a_0}{\chi + \lambda_2 + 2\bar{\gamma}\tau}$. Hence, for $\lambda_{sup} = \lambda_2$, we obtain, $\rho_{SM}(g^*) = 1 - \frac{2\bar{\gamma}\tau\lambda_2\left(g^*\right)^2}{\left(1 - \frac{1}{\delta}\right)^2 - \left(a_0 - \lambda_2g^*\right)^2}$.

\section{Acknowledgments}
This work is supported by National Science Foundation ECCS grants 1002053, 1150405, and CNS grant 1329915 (to U.V.).

\section{Contribution of Authors}
U.V formulated the problem and A. V proved the main results.  A.D. and U.V. wrote the main text of the manuscript. A.D. ran the simulations and prepared all the figures. Both A.D. and U.V. reviewed the manuscript.

\section{Additional Information}
{\bf Competing financial interests:} The authors declare no competing financial interests.


\begin{thebibliography}{0}
\expandafter\ifx\csname natexlab\endcsname\relax\def\natexlab#1{#1}\fi
\expandafter\ifx\csname bibnamefont\endcsname\relax
  \def\bibnamefont#1{#1}\fi
\expandafter\ifx\csname bibfnamefont\endcsname\relax
  \def\bibfnamefont#1{#1}\fi
\expandafter\ifx\csname citenamefont\endcsname\relax
  \def\citenamefont#1{#1}\fi
\expandafter\ifx\csname url\endcsname\relax
  \def\url#1{\texttt{#1}}\fi
\expandafter\ifx\csname urlprefix\endcsname\relax\def\urlprefix{URL }\fi
\providecommand{\bibinfo}[2]{#2}
\providecommand{\eprint}[2][]{\url{#2}}

\end{thebibliography}


\begin{thebibliography}{99}

\bibitem{Strogatz2}
Strogatz,~S.~H. \& Stewart,~I.
\newblock Coupled oscillators and biological synchronization.
\newblock {\em ,Sci. Am.}, {\bf 269}, 102--109(1993).

\bibitem{kuramoto_bonilla1}
Acebr\'on,~J.~A., Bonilla,~L.~L., P\'erez~Vicente,C.~J., Ritort,~F. \& Spigler,R.
\newblock The kuramoto model: A simple paradigm for synchronization phenomena.
\newblock {\em Rev. Mod. Phys.}, {\bf 77}, 137--185(2005).

\bibitem{genetic_clock_hasty}
Danino,~T., Mondragon-Palomino,~O., Tsimring,~L., \& Hasty,~J.
\newblock A synchronized quorum of genetic clocks.
\newblock {\em Nature}, {\bf 463}, 326--330(2010).

\bibitem{PNAS_Dorfler}
D\"orfler,~F., Chertkov,~M., \& Bullo,~F.
\newblock Synchronization in complex oscillator networks and smart grids.
\newblock {\em Proc. Natl. Acad. Sci. U.S.A.}, {\bf 110}, 2005--2010(2013).

\bibitem{Chaos_Rohden}
Rohden,~M., Sorge,~A., Witthaut,~D., \& Timme,~M.
\newblock Impact of network topology on synchrony of oscillatory power grids.
\newblock {\em Chaos}, {\bf 24}, 013123(2014),
{\em http://dx.doi.org/10.1063/1.4865895}.

\bibitem{Chaos_Stout}
Stout,~J., Whiteway,~M., Ott,~E., Girvan,~M., \& Antonsen,~T.~M.
\newblock Local synchronization in complex networks of coupled oscillators.
\newblock {\em Chaos}, {\bf 21}, 025109(2011),
{\em http://dx.doi.org/10.1063/1.3581168}.

\bibitem{Erdos_Renyi}
Erd\"{o}s,~P., \& R\'{e}nyi,~A.
\newblock {On random graphs, I}.
\newblock {\em Publ. Math-Debrecen}, {\bf 6}, 290--297(1959).

\bibitem{Nature_watts}
Watts,~D.~J., \& Strogatz,~S.~H.
\newblock {Collective dynamics of small-world networks.}
\newblock {\em Nature}, {\bf 393}, 409--10(1998).

\bibitem{PNAS_Amaral}
Amaral,~L.~A.~N., Scala,~A., Barth\'el\'emy,~M., \& Stanley,~H.~E.
\newblock Classes of small-world networks.
\newblock {\em Proc. Natl. Acad. Sci. U.S.A.}, {\bf 97}, 11149--11152(2000).

\bibitem{Science_Barabasi}
Barab\'asi,~A.~L., \& Albert,~R.
\newblock Emergence of scaling in random networks.
\newblock {\em Science}, {\bf 286}, 509--512(1999).

\bibitem{Pecora_Nature}
Pecora,~L.~M., Sorrentino,~F., Hagerstrom,~A.~M., Murphy,~T.~E., \& Roy,~R.
\newblock Cluster synchronization and isolated desynchronization in complex networks with symmetries.
\newblock {\em Nat. Comm.}, {\bf 5}, (2014),
{\em doi:10.1038/ncomms5079}.

\bibitem{Arndt_Nature}
Becks,~L., \& Arndt,~H.
\newblock Different types of synchrony in chaotic and cyclic communities.
\newblock {\em Nat. Comm.}, {\bf 4}, 1359-1367(2013).

\bibitem{pecora_carroll_master}
Pecora,~L.~M., \& Carrol,~T.~L.
\newblock Master stability functions for synchronized coupled systems.
\newblock {\em Phys. Rev. Lett.}, {\bf 80}, 2109--2112(1998).

\bibitem{pecora_barahona}
Barahona,~M. \& Pecora,~L.~M.
\newblock Synchronization in Small-World systems.
\newblock {\em Phys. Rev. Lett.}, {\bf 89}, 0112023(2002).

\bibitem{rangarajan_ding}
Rangarajan,~G., \& Ding,~M.
\newblock Stability of synchronized chaos in coupled dynamical systems.
\newblock {\em Phys. Lett. A}, {\bf 296}, 159--187(2002).

\bibitem{belykh_connection}
Belykh,~V., Belykh,~I., \& Hasler,~M.
\newblock Connection graph stability method for synchronized coupled chaotic systems.
\newblock {\em Physica D}, {\bf 195}, 159--187(2004).

\bibitem{PNAS_Mezic}
Mezic,~I.
\newblock On the dynamics of molecular conformation.
\newblock {\em Proc. Natl. Acad. Sci. U.S.A.}, {\bf 103}, 7542--7547(2006).

\bibitem{PNAS_Nishikawa}
Nishikawa,~T., \& Motter,~A.~E.
\newblock Network synchronization landscape reveals compensatory structures,
  quantization, and the positive effect of negative interactions.
\newblock {\em Proc. Natl. Acad. Sci. U.S.A.}, {\bf 107}, 10342--10347(2010).

\bibitem{synchronization_scalefree}
Wang,~X., \& Chen,~G.
\newblock Synchronization in scale-free dynamical networks: Robustness and fragility.
\newblock {\em IEEE Trans. Circuits Syst. I, Fundam. Theory Appl.}, {\bf 49}, 54--62(2002).

\bibitem{porfiri_master}
Porfiri,~M.
\newblock A master stability function for stochastically coupled chaotic maps.
\newblock {\em Europhys. Ltt.}, {\bf 6}, 40014(2011).

\bibitem{amit_synchronization}
Diwadkar,~A., \& Vaidya,~U.
\newblock Robust synchronization in network systems with link failure uncertainty.
\newblock In {\em {Proc. IEEE Decis. Contr. P.}}, 6325--6330(2010).

\bibitem{belykh_blinking}
Belykh,~I., Belykh,~V., \& Hasler,~M.
\newblock Blinking models and synchronization in small-world networks with a time-varying coupling.
\newblock {\em Physica D}, {\bf 195}, 188--206(2004).

\bibitem{Belykh_SIAM_I}
Hasler,~M. \& Belykh,~V. \& Belykh,~I.
\newblock Dynamics of Stochastically Blinking Systems. Part I: Finite Time Properties.
\newblock In {\em {SIAM J. Appl. Dyn. Syst.}}, {\bf 12}, 1007--1030(2013).

\bibitem{Belykh_SIAM_II}
Hasler,~M. \& Belykh,~V. \& Belykh,~I.
\newblock Dynamics of Stochastically Blinking Systems. Part II: Asymptotic Properties.
\newblock In {\em {SIAM J. Appl. Dyn. Syst.}}, {\bf 12}, 1031--1084(2013).

\bibitem{Jost_SIAM}
Lu,~W \& Atay,~F.~M. \& Jost,~J.
\newblock Synchronization of discrete-time dynamical networks with time-varying coupling
\newblock In {\em {SIAM J. Math. Anal.}}, {\bf 39}, 1231--1259(2007).

\bibitem{Jost_Atay}
Lu,~W \& Atay,~F.~M. \& Jost,~J.
\newblock Chaos synchronization in networks of coupled maps with time-varying topologies
\newblock In {\em {SIAM J. Math. Anal.}}, {\bf 63}, 399--406(200).

\bibitem{PRE_Molino}
Garcia~del Molino,~L.~C., Pakdaman,~K., Touboul,~J. \& Wainrib,~G.
\newblock Synchronization in random balanced networks.
\newblock {\em Phys. Rev. E}, {\bf 88}, 042824(2013).

\bibitem{PRE_Sinha}
Sinha,~S., \& Sinha,~S.
\newblock Robust emergent activity in dynamical networks.
\newblock {\em Phys. Rev. E}, {\bf 74}, 066117(2006).

\bibitem{PRE_Kocarev}
Kocarev,~L., Parlitz,~U., \& Brown,~R.
\newblock Robust synchronization of chaotic systems.
\newblock {\em Phys. Rev. E}, {\bf 61}, 3716--3720(2000).

\bibitem{PRE_Wang}
Wang,~Z., Fan,~H., \& Aihara,~K.
\newblock Three synaptic components contributing to robust network synchronization.
\newblock {\em Phys. Rev. E}, {\bf 83}, 051905(2011).

\bibitem{amit_erasure_observation_journal}
Diwadkar,~A., \& Vaidya,~U.
\newblock Limitation on nonlinear observation over erasure channel.
\newblock {\em {IEEE Trans. Autom. Control}}, {\bf 58}, 454--459(2013).

\bibitem{amit_ltv_journal}
Diwadkar,~A., \& Vaidya,~U.
\newblock Stabilization of linear time varying systems over uncertain channels.
\newblock {\em {Int. J. Robust Nonlin.}}, {\bf 24}, 1205-1220(2014).

\bibitem{Vaidya_erasure_SCL}
Vaidya,~U., \& Elia,~N.
\newblock Limitation on nonlinear stabilization over packet-drop channels:
  Scalar case.
\newblock {\em {Syst. Control Lett.}}, {\bf 61}, 959--966(2012).

\bibitem{Vaidya_erasure_stablization}
Vaidya,~U., \& Elia,~N.
\newblock Limitation on nonlinear stabilization over erasure channel.
\newblock In {\em {Proc. IEEE Decis. Contr. P.}}, 7551--7556(2010).

\bibitem{Hasminskii_book}
Has\'minski\u{i},~R.~Z.
\newblock {\em Stability of differential equations}.
\newblock Sijthoff \& Noordhoff, {Germantown ,MD}, (1980).

\bibitem{Wang_IEEE_TNN}
Wang,~Z., Wang,~Y., \& Liu,~Y.
\newblock Global synchronization for discrete-time stochastic complex networks with randomly occurred nonlinearities and mixed time delays.
\newblock {\em IEEE Trans. Neural Netw.}, {\bf 21}, 11--25(2010).

\bibitem{Astrom_book}
Astrom,~K.~J., \& Murray,~R.~M.
\newblock {\em Feedback Systems: An Introduction for Scientists and Engineers}.
\newblock Princeton University Press, {Princeton, NJ}, (2008).

\bibitem{scl04}
Elia,~N,
\newblock Remote Stabilization over Fading Channels.
\newblock {\em Syst. Control Lett.}, {\bf 54}, 237-249(2005)

\bibitem{Haddad_book}
Haddad,~W., \& Chellaboina,~V.~S.
\newblock {\em Nonlinear dynamical systems and control: A Lyapunov-based approach}.
\newblock Princeton University Press, {Princeton, NJ}, (2008).

\bibitem{haddad_bernstein_DT}
Haddad,~W., \& Bernstein,~D.
\newblock Explicit construction of quadratic {L}yapunov functions for the small gain theorem, positivity, circle and {P}opov theorems and their application to robust stability. {P}art {II}: {D}iscrete-time theory.
\newblock {\em Int. J. Robust Nonlin.}, {\bf 4}, 249--265(1994).

\bibitem{lancaster}
Lancaster,~P., \& Rodman,~L.
\newblock {\em Algebraic Riccati Equations}.
\newblock Oxford Science Publications, Oxford, (1995).

\bibitem{Merris1994}
Merris,~R.
\newblock Laplacian matrices of graphs: a survey.
\newblock {\em Linear Algebra Appl.}, {\bf 197-198}, 143--176(1994).

\bibitem{Boyd_book}
Boyd,~S., \& Vandenberghe,~L.
\newblock {\em Convex Optimization}.
\newblock Cambridge University Press, {Cambridge, UK}, (2003).

\end{thebibliography}
\end{document}


\title{Limitations and tradeoff in synchronization of large-scale stochastic networks}

\author{Amit Diwadkar and Umesh Vaidya}
\affiliation{Electrical and Computer Engineering, Iowa State University \\ Coover Hall, Ames, IA, U.S.A. 50011}

\maketitle

\section{Supplementary Information}
We consider the problem of synchronization in large-scale nonlinear network systems, where the scalar dynamics of the individual subsystem is assumed,
\begin{equation} \label{component_dynamics}
x^k _{t+1} = a x^k_t - {\phi}  (x^k_t) + v^k_t ~ ~k=1,\ldots, N,
\end{equation}
where $x^k \in \mathbb{R}$ are the states of the $k^{th}$ subsystem, $a > 0$, and $v^k \in \mathbb{R}$ is an independent identically distributed (i.i.d.) additive noise process with zero mean (i.e., $E[v^k_t] = 0$) and variance $E[(v^k_t)^2] = \omega^2$.  Subscript $t$ denotes the index of the discrete time-step throughout the paper. The function, $\phi \colon \mathbb{R} \to \mathbb{R}$, is a monotonic, globally Lipschitz function with $\phi(0) = 0$ and Lipschitz constant, $\frac{2}{\delta}$. If the coupling gain is $g > 0$, the individual agent dynamics of the coupled subsystem is given by,
\begin{align}
x^k_{t+1} &= a x^k_t - \phi\left(x^k_t\right) + g\sum_{e_{kj}\in \Ed_D}\mu_{kj} (x^j_t - x^k_t) +g\sum_{e_{kj} \in \Ed\Ed_U}(\mu_{kj}+\xi_{kj}) (x^j_t - x^k_t) + v^k_t,\\
&= a x^k_t - \phi\left(x^k_t\right) + g\sum_{e_{kj}\in \Ed}\mu_{kj} (x^j_t - x^k_t) +g\sum_{e_{kj} \in \Ed_U}\xi_{kj} (x^j_t - x^k_t) + v^k_t,\\
&= a x^k_t - \phi\left(x^k_t\right) - g\left(\sum_{e_{kj}\in \Ed}\mu_{kj}\right)x^k_t + g\sum_{e_{kj}\in \Ed}\mu_{kj} x^j_t - g\left(\sum_{e_{kj} \in \Ed_U}\xi_{kj}\right)x^k_t +g\sum_{e_{kj} \in \Ed_U}\xi_{kj}x^j_t + v^k_t.
\end{align}
We denote the nominal graph Laplacian by ${\cal L} := \left[ l(ij)\right] \in \R^{N \times N}$ and uncertain graph Laplacian by ${\cal L}_R := \left[ l_R(ij)\right] \in \R^{N \times N}$, where 
\begin{align} 
l(ij) = \left\{\begin{array}{cc} 
-\mu_{ij}, & \text{if $i\neq j$, and, $e_{ij} \in \Ed$}\\  
\sum_{e_{ij}\in \Ed}\mu_{ij}, & \text{if $i=j$}
\end{array}\right.,\qquad \text{and}, \qquad
l_R(ij) = \left\{\begin{array}{cc} 
-\xi_{ij}, & \text{if $i\neq j$, and, $e_{ij} \in \Ed_U$}\\  
\sum_{e_{ij}\in \Ed_U}\xi_{ij}, & \text{if $i=j$}
\end{array}\right. .
\end{align}
We combine the individual systems to create the network system ($\tilde{x}_t$) written as,
\begin{align}
\tilde{x}_{t+1} = \left(aI_N - g({\cal L}+{\cal L}_R)\right) \tilde{x}_t - \tilde{\phi} \left (\tilde x_t\right) + \tilde{v}_t, \label{coupled_dynamics}
\end{align}
where $I_N$ is the $N\times N$ identity matrix, $\tilde x_t  = [ x^1_t\; \cdots\; x^N_t ]^{\top}$, and $\tilde {\phi}(\tilde x_t) = [ {\phi}^1_t(x^1_t) \; \cdots \; {\phi}^N_t(x^N_t) ]^{\top}$.

Given the stochastic nature of the network system, we propose the following definition for mean square exponential synchronization \cite{Hasminskii_book}.
\begin{definition}[Mean Square Synchronization]\label{mse_sync_def}
The network system \eqref{coupled_dynamics} is said to be mean square synchronizing (MSS), if there exist  positive constants, $\beta < 1$, $\bar{K}(\tilde e_0) <\infty$, and $L <\infty$, such that,
\begin{eqnarray}
& E_{\Xi} \parallel x^k_{t} - x^j_{t}\parallel^2 \leq \bar K(\tilde e_0) {\beta}^t \parallel x^k_{0} - x^j_{0}\ \parallel^2 + L\omega^2,
\label{mse_sync_eqn}
\end{eqnarray}
 $\forall k,j \in [1,N]$, where $\tilde e_0$ is a function of difference $\parallel x_0^i-x_0^\ell\parallel^2$ for $i,\ell \in [1,N]$ and $\bar K(0)=K$ for some constant $K$.
\end{definition}
\begin{remark} In the absence of additive noise, $\tilde v_t$, in system Eq. \eqref{coupled_dynamics}, the term  $L \omega^2$ in Eq. \eqref{mse_sync_eqn} vanishes and Definition \ref{mse_sync_def} then reduces to mean square exponential (MSE) synchronization \cite{Wang_IEEE_TNN}.
\end{remark}

\subsection{Mean Square Synchronization Result}

Since the subsystems are identical, the synchronization manifold is spanned by the vector, ${\mathds 1}=[1,\ldots,1]^{\top}$. The dynamics on the synchronization manifold are decoupled from the dynamics off the manifold and are essentially described by the dynamics of the individual system, which could be stable, oscillatory, or complex in nature. We now apply a change of coordinates to decompose the system dynamics on and off the synchronization manifold.
Let ${\cal L} = V\Lambda V^{\top}$, where $V$ is an orthonormal set of vectors given by $V = \left[\frac{\mathds{1}}{\sqrt{N}}\;\; U\right]$, and $U$ is a set of orthonormal vectors also orthonormal to $\mathds{1}$. Furthermore, $\Lambda = \text{diag}\{\lambda_1,\cdots,\lambda_N\}$, where $0 = \lambda_1 < \lambda_2 \leq \cdots \leq \lambda_N$ are the eigenvalues of ${\cal L}$. Let $\tilde z_t = V^{\top}\tilde x_t$ and $\tilde w_t = V^{\top}\tilde v_t$. Multiplying (\ref{coupled_dynamics}) from the left side by $V^{\top}$, we obtain
\begin{align}
\tilde{z}_{t+1} = \left( aI_N - g\left( V^{\top}({\cal L}+{\cal L}_R)V \right)\right ) \tilde{z}_t - \tilde{\psi} \left (\tilde z_t\right ) + \tilde w_t, \label{coupled_dyn_trans}
\end{align}
where $\tilde{\psi}(\tilde z_t) = V^{\top}\tilde{\phi}\left(\tilde x_t\right)$. We can now write $\tilde z_t = \left[\begin{array}{cc}\bar{x}_t^{\top} & \hat{z}_t^{\top} \end{array}\right]^{\top}$, $\tilde{\psi}(\tilde z_t) := \left[\begin{array}{cc} \bar{\phi}_t^{\top} & \hat{\psi}_t^{\top} \end{array}\right]^{\top}$,
and 
$
\tilde{w}_t := \left[\begin{array}{cc}
\bar{v}_t^{\top} & \hat{w}_t^{\top}
\end{array}\right]^{\top}
$, where
\begin{align}
\bar{x}_t &:= \frac{\mathds{1}^{\top}}{\sqrt{N}}\tilde{x}_t = \frac{1}{\sqrt{N}}\sum_{k = 1}^N x^k_t,\quad \hat{z}_t := U^{\top}\tilde{x}_t \\
\bar{\phi}_t &:= \frac{\mathds{1}^{\top}}{\sqrt{N}}\tilde{\phi}\left(\tilde z_t\right) = \frac{1}{\sqrt{N}}\sum_{k = 1}^N \phi(x^k_t),\quad \hat{\psi}_t := U^{\top}\tilde{\phi}\left(\tilde x_t\right).
\end{align}
Furthermore, we have
\begin{align}
E_{\tilde v}[\bar{v}_t^2] = \sqrt{N}\omega^2,\; E_{\tilde v}[\hat{w}_t\hat{w}_t^{\top}] = U^{\top}E_{\tilde v}[\tilde{v}_t\tilde{v}_t^{\top}]U = \omega^2I_{N-1}.
\end{align}
From \eqref{coupled_dyn_trans}, we obtain
\begin{align}
\bar{x}_{t+1} &=  a \bar{x}_t - \bar{\phi} \left (\bar x_t\right ) + \bar{v}_t\nonumber\\
\hat{z}_{t+1} &= \left( aI_{N-1} -  g\left( \hat{\Lambda} + U^{\top}{\cal L}_R U \right)\right ) \hat{z}_t - \hat{\psi}_t + \hat{w}_t, \label{coupled_dyn_reduced}
\end{align}
where $\hat{\Lambda} = \text{diag}\{\lambda_2,\cdots,\lambda_N\}$. For the synchronization of system (\ref{coupled_dynamics}), we only need to demonstrate the mean square stability to the origin of the $\hat{z}$ dynamics, as given in (\ref{coupled_dyn_reduced}). This feature is exploited to derive the sufficiency condition for mean square synchronization of the coupled system, as shown in the following lemma.
\begin{lemma}\label{MSEsync_mss}
The system described by Eq. \eqref{coupled_dynamics} is MSS, as given by Definition 1, if there exists $L > 0$, $K > 0$, and $0 < \beta < 1$, such that
\begin{align}
E_{\Xi} \left[ \parallel \hat z_{t} \parallel^2 \right] &\le   K  {\beta}^t \parallel \hat z_{0}  \parallel^2 + L\omega^2.
\end{align}
\end{lemma}
\begin{proof} To prove this result, we show the second moment of $\hat z_t$ dynamics is equivalent to the mean square error dynamics for each pair of systems. Then, we apply the stability results to the error dynamics to complete the proof. Consider Eq. \eqref{coupled_dyn_reduced}. We have
\begin{align}
\parallel \hat z_t \parallel^2 &= \hat z_t^{\top} \hat z_t = \tilde x_t^{\top} \left(UU^{\top} \otimes I_n\right)\tilde x_t. \label{diff_eq1}
\end{align}
Then, we have $UU^{\top} = VV^{\top} - \frac{\bf{1}}{\sqrt{N}}\frac{\bf{1}^{\top}}{\sqrt{N}} = I_N - \frac{1}{N}\bf{1}\bf{1}^{\top}$. Substituting in \eqref{diff_eq1}, we obtain
\begin{align}
\parallel \hat z_t \parallel^2 = \frac{1}{2N}\sum_{i=1}^N\sum_{j\neq i,j=1}^N \left(x_t^i - x_t^j\right)^{\top}\left(x_t^i - x_t^j\right).
\end{align}
Now, suppose there exist $L > 0$, $K > 0$, and $0 < \beta < 1$, such that
\begin{align}
E_{\Xi} \left[\parallel \hat z_{t} \parallel^2 \right] &\le   K  {\beta}^t \parallel \hat z_{0}  \parallel^2 + L\omega^2.
\end{align}
This can be rewritten as
\begin{align}
&E_{\Xi} \left[\sum_{k=1}^N\sum_{j\neq k, j=1}^N \parallel x^k_t - x^j_t \parallel^2 \right] \le  K  {\beta}^t  \sum_{k=1}^N\sum_{j\neq k,j=1}^N \parallel x^k_0 - x^j_0  \parallel^2 + L\omega^2.
\end{align}
This implies
\begin{align}\label{sum_error}
&\sum_{k=1}^N\sum_{j\neq k, j=1}^N E_{\Xi}\left[\parallel x^k_t -x^j_t \parallel^2 \right] \le  K  {\beta}^t  \sum_{k=1}^N\sum_{j\neq k,j=1}^N \parallel x^k_0 - x^j_0  \parallel^2 + L\omega^2.
\end{align}
Thus, from \eqref{sum_error} we obtain for all systems, $S_k$ and $S_l$,
\begin{align}
E_{\Xi} \parallel x^k_t - x^l_t \parallel^2 \le  \bar K (\tilde e_0) \beta^t \parallel x^k_0 - x^l_0 \parallel^2 + L\omega^2,
\end{align}
where $\bar K(\tilde{e}_0) := K \left ( 1 + \frac{\sum_{i=1, i\neq k}^N\sum_{j=1, j\neq i}^N \parallel x^i_0 - x^j_0 \parallel^2}{\parallel x^k_0 - x^l_0 \parallel^2}  \right )$. Hence, the proof.
\end{proof}

We will now provide a slightly modified Eq. \eqref{coupled_dyn_reduced}. We know ${\cal L}_R = \sum_{e_{ij}\in \Ed_U}\xi_{ij}\ell_{ij}\ell_{ij}^{\top}$, where $\ell_{ij} \in \mathbb{R}$ has values $1$ and $-1$ in $i^{th}$ and $j^{th}$ entries, respectively, the remaining values are zeros. Thus, $\ell_{ij}^{\top}\ell_{ij} = 2$ for all $e_{ij} \in \Ed$. Hence, if $\hat{\ell}_{ij} = U^{\top}\ell_{ij}$, we calculate
\begin{align}
U^{\top}{\cal L}_R U_d = \sum_{e_{ij}\in \Ed_U}\xi_{ij}U^{\top}\ell_{ij}\ell_{ij}^{\top}U = \sum_{e_{ij}\in \Ed_U}\xi_{ij}\hat{\ell}_{ij}\hat{\ell}_{ij}^{\top},
\end{align}
where $\hat{\ell}_{ij}^{\top}\hat{\ell}_{ij} = 2$. Thus, we can write Eq. \eqref{coupled_dyn_reduced} as
\begin{align}
\hat z_{t+1} = \left(aI_{N-1} - g\hat{\Lambda} - \sum_{e_{ij}\in \Ed_U}\xi_{ij}\hat{\ell}_{ij}\hat{\ell}_{ij}^{\top}\right)\hat{z}_t - \hat{\psi}(\hat{z}_t) + \hat{w}_t.
\label{sync_dynamics}
\end{align}
In Lemma \ref{MSEsync_mss}, we prove mean square exponential stability of \eqref{sync_dynamics} guarantees the mean square synchronization of the coupled network of Lure systems, as given by \eqref{coupled_dynamics}. We will now utilize this to provide the sufficiency condition for mean square stabilization of Lure systems interacting over a network.

\begin{theorem}\label{1D_suff_original}
The network system in Eq. \eqref{coupled_dynamics} is mean square synchronizing, if there exists a positive constant, $p$, that satisfies $\delta > p$,
\begin{align}\label{sync_ric}
\left( p - \frac{1}{\delta} \right)\left(\frac{1}{p} - \frac{1}{\delta} \right) > \alpha_0^2,
\end{align}
where $\alpha_0^2 = (a_0 - \lambda_{sup} g)^2 + 2\bar{\gamma}\tau\lambda_{sup}g^2$, $a_0 = a - \frac{1}{\delta}$, and $\lambda_{sup} = \underset{\lambda \in \{\lambda_2,\lambda_N\}}{argmax}\bigg|\lambda+\bar{\gamma}\tau-\frac{a_0}{g}\bigg|$.
Furthermore, $\tau := \frac{\lambda_{N_U}}{\lambda_{N_U}+\lambda_{2_D}}$, where $\lambda_{N_U}$ is the maximum eigenvalue of ${\cal L}_U$, and $\lambda_{2_D}$ is the second smallest eigenvalue of ${\cal L}_D$.
\end{theorem}

\begin{proof}
We first construct an appropriate Lyapunov function, $V(\hat{z}_t) = \hat{z}_t^{\top}P\hat{z}_t$, that guarantees mean square stability. From \eqref{sync_dynamics}, defining $\Delta V :=E_{\Xi_t,v_t}\left[ V(\hat{z}_{t+1}) - V(\hat{z}_t) \right]$, we obtain,
\begin{align}
\Delta V &= E_{\Xi_t}\left[\hat{z}_t^{\top}\left(A(\Xi_t)^{\top}PA(\Xi_t) - P\right)\hat{z}_t -\hat{z}_t^{\top}A(\Xi_t)^{\top}P\hat{\psi}_t - \hat{\psi}_t^{\top}PA(\Xi_t)\hat{z}_t + \hat{\psi}_t^{\top}P\hat{\psi}_t\right]
 + E_{v_t}[\hat{w}_t^{\top}P\hat{w}_t].
\end{align}
Now suppose for some $R_P > 0$, $P$ satisfies
\begin{align}\label{unc_lure_p1}
P =& E_{\Xi} \left [ A(\Xi)^{\top}PA(\Xi)\right]  + R_P + E_{\Xi} \left [\left(A(\Xi)^{\top}P-I_{N-1}\right)(\delta I_{N-1}-P)^{-1}\left(P A(\Xi)-I_{N-1}\right)\right ].
\end{align}
Using \eqref{unc_lure_p1} and algebraic manipulations as given in \cite{haddad_bernstein_DT}, we can rewrite $\Delta V$ as
\begin{align}
\Delta V &= - \hat{z}^{\top}_t R_P \hat{z}_t - E_{\Xi_t}\left[\eta^{\top}_t \eta_t \right] - 2 \hat\psi^{\top}_t\left(\hat{z}_t - \frac{\delta}{2}\hat\psi_t\right) + \textrm{trace}(P E_{v_t}[\hat{w}_t\hat{w}_t^{\top}]),
\end{align}
where $\eta_t(\Xi(t))$ be given by $\eta_t(\Xi(t)) = W^{-\frac{1}{2}}\left(PA(\Xi) - I_{N-1}\right)\hat{z}_t - W^{\frac{1}{2}}\hat\psi_t$ and $W := (\delta I_{N-1} - P)$. 
Since $\phi(\cdot)$ is monotonic and globally Lipschitz with constant $\frac{2}{\delta}$, we know $\left(\phi\left(x^k_t\right)-\phi\left(x^l_t\right)\right)^{\top}\left(\frac{2}{\delta}\left(x_t^k - x_t^l\right) - \left(\phi\left(x^k_t\right)-\phi\left(x^l_t\right)\right)\right) > 0$. 
This gives 
\begin{align}\label{nonlin_ineq}
\hat{\psi}_t^{\top}\left(\hat{z}_t - \frac{\delta}{2}\hat{\psi}_t\right) &= \tilde{\phi}\left(\tilde{x}_t\right)^{\top}UU^{\top}\left(\frac{2}{\delta}\tilde{x}_t - \tilde{\phi}\left(\tilde{x}_t\right)\right) = \tilde{\phi}\left(\tilde{x}_t\right)^{\top}\left( I - \frac{1}{N}\mathds{1}\mathds{1}^{\top}\right)\left(\frac{2}{\delta}\tilde{x}_t - \tilde{\phi}\left(\tilde{x}_t\right)\right) \\
&= \frac{1}{2N}\sum_{k=1}^N\sum_{l=1,l\neq k}^N \left(\phi\left(x^k_t\right)-\phi\left(x^l_t\right)\right)^{\top}\left(\frac{2}{\delta}\left(x_t^k - x_t^l\right) - \left(\phi\left(x^k_t\right)-\phi\left(x^l_t\right)\right)\right) > 0.
\end{align}
Using \eqref{nonlin_ineq} and writing $\rho = \textrm{trace}(P)$, we obtain,
\begin{align}\label{1step_lyap2} 
E_{\Xi_t,v_t} \left  [ V(\hat{z}_{t+1})) - V(\hat{z}_t) \right ] < - \hat{z}^{\top}_t R_P \hat{z}_t + \rho \omega^2.
\end{align} 
Hence, \eqref{unc_lure_p1} is sufficient for MSS of \eqref{component_dynamics} from condition for Mean Square Stability of the Reduced System. Furthermore, the equation in (\ref{unc_lure_p1}) can be rewritten using ~\cite{lancaster} (Proposition 12.1,1) as
\begin{align} \label{unc_lure_p}
P =& E_{\Xi} \left [ A_0(\Xi)^{\top}PA_0(\Xi)\right]  + R_P + \frac{1}{\delta}I_{N-1}+E_{\Xi} \left [A_0(\Xi)^{\top}P(\delta I_{N-1}-P)^{-1}P A_0(\Xi)\right ],
\end{align}
where $A_0(\Xi) = a_0I_{N-1} -g\hat{\Lambda} - gU^{\top}{\cal L}_R U$ and $a_0 = a - \frac{1}{\delta}$. We observe this condition requires us to find a symmetric Lyapunov function matrix, $P$, of order $\frac{N(N-1)}{2}$. We now reduce the order of computation by using network properties. For this, consider $P = p I_{N-1}$, where $p < \delta$ is a positive scalar. This gives us $\delta I_{N-1} > P$. Using this and \eqref{sync_dynamics}, we rewrite the condition in \eqref{unc_lure_p} as follows
\begin{align} \label{ric1}
pI_{N-1} &>  p(a_0I_{N-1} - g\hat{\Lambda})^{\top}(a_0I_{N-1}- g\hat{\Lambda})  + \frac{p^2}{\delta - p}(a_0I_{N-1}-g\hat{\Lambda})^{\top}(a_0I_{N-1}-g\hat{\Lambda}) + \frac{1}{\delta}I_{N-1}\nonumber \\
&\quad+ pg^2\underset{e_{ij \in \Ed_U}}{\sum} \sigma_{ij}^2\hat{\ell}_{ij}\hat{\ell}_{ij}^{\top}\hat{\ell}_{ij}\hat{\ell}_{ij}^{\top} +\frac{p^2g^2}{\delta - p}\underset{e_{ij \in \Ed_U}}{\sum} \sigma_{ij}^2\hat{\ell}_{ij}\hat{\ell}_{ij}^{\top}\hat{\ell}_{ij}\hat{\ell}_{ij}^{\top}.
\end{align}
We know $\hat{\ell}_{ij}^{\top}\hat{\ell}_{ij} = \ell_{ij}^{\top}U_dU_d^{\top}\ell_{ij} = \ell_{ij}^{\top}\ell_{ij} = 2$ and
\begin{align}\label{co_dis}
\underset{e_{ij \in \Ed_U}}{\sum}\sigma_{ij}^2\hat{\ell}_{ij}\hat{\ell}_{ij}^{\top}\hat{\ell}_{ij}\hat{\ell}_{ij}^{\top}\leq 2\bar{\gamma}\underset{e_{ij \in \Ed_U}}{\sum}\mu_{ij}\hat{\ell}_{ij}\hat{\ell}_{ij}^{\top} = 2\bar{\gamma}U^{\top}{\cal L}_U U.
\end{align}
Now, suppose 
${\cal L}_U\leq \tau {\cal L} = \tau \left({\cal L}_D + {\cal L}_U\right)$. We then obtain $\tau \geq \frac{\lambda_{N_U}}{\lambda_{N_U}+\lambda_{2_D}}$. We choose $\tau$ as
\begin{align}\label{tau_value}
\tau = \frac{\lambda_{N_U}}{\lambda_{N_U}+\lambda_{2_D}}.
\end{align}
Hence, we have ${\cal L}_U\leq \tau \left({\cal L}_D + {\cal L}_U\right) = \tau {\cal L}$. Then, bounding ${\cal L}_U$ in \eqref{co_dis} using \eqref{tau_value}, we obtain,
\begin{align}\label{co_dis2}
\underset{e_{ij \in \Ed_D}}{\sum}\sigma_{ij}^2\hat{\ell}_{ij}\hat{\ell}_{ij}^{\top}\hat{\ell}_{ij}\hat{\ell}_{ij}^{\top} \leq 2\bar{\gamma}\tau U^{\top} {\cal L} U = 2\bar{\gamma}\tau \hat{\Lambda}.
\end{align}
Substituting \eqref{co_dis2} into \eqref{ric1}, a sufficient condition for inequality \eqref{ric1} to hold is given by
\begin{align} \label{ric4}
 pI_{N-1} > \left(p+\frac{p^2}{\delta - p}\right)(a_0I_{N-1} - g\hat{\Lambda})^{\top}(a_0I_{N-1}- g\hat{\Lambda} ) + \left(p+\frac{p^2}{\delta - p}\right)2\bar{\gamma}\tau g^2\hat{\Lambda} + \frac{1}{\delta}I_{N-1}.
\end{align}
Equation \eqref{ric4} is a diagonal matrix equation that provides a sufficient condition for MSS as
\begin{align} \label{ric5}
p >  \left(p+\frac{p^2}{\delta - p}\right)\left((a_0 - g\lambda_j)^2 + 2\bar{\gamma}\tau g^2\lambda_j\right) + \frac{1}{\delta},
\end{align}
for all eigenvalues $\lambda_j$ of $\hat{\Lambda}$. This is simplified as
\begin{align}\label{suff_main1}
\left(p - \frac{1}{\delta}\right)\left(\frac{1}{p} - \frac{1}{\delta}\right) > \alpha_0^2,
\end{align}
where $\delta > p > 0$ and $\alpha_0^2 = (a_0 - g\lambda)^2 + 2\bar{\gamma}\tau \lambda g^2$ for all $\lambda \in \{\lambda_2,\ldots,\lambda_N\}$ are eigenvalues of the nominal graph Laplacian. Now, for each of these conditions to hold true, we must satisfy condition \eqref{suff_main1} for the minimum value of $\alpha_o^2$ with respect to all possible $\lambda$. Now, $\lambda^*$ provides minimum values for $\alpha_0^2$ found by setting $\frac{d\alpha_0^2}{d\lambda}\Big|_{\lambda^*} = 0$, giving us 
\begin{align}
\lambda^* = \frac{a_0}{g} - \bar{\gamma}\tau.
\end{align}
Using $\lambda^*$, we conclude, if \eqref{suff_main1} is to be satisfied for all $\lambda \in \{\lambda_2\ldots,\lambda_N\}$, it must satisfy \eqref{suff_main1} for the farthest such $\lambda$ from $\lambda^*$. Since eigenvalues of a graph Laplacian are positive and monotonic non-decreasing, all we need is to satisfy \eqref{suff_main1} for $\lambda_{sup}$, where $\lambda_{sup} = \underset{\lambda \in \{\lambda_2,\lambda_N\}}{\text{argmax}}|\lambda - \lambda^*|$.
\end{proof}

In the following discussion we will provide a system theoretic interpretation to the proposed definition of mean square synchronization margin. For completion, we restate the definition for mean square synchronization margin.
\begin{definition}[Mean Square Synchronization Margin] \label{def_sm}The margin for synchronization for network system  \eqref{coupled_dynamics} is defined as
\begin{align}\label{sync_margin}
\rho_{SM} := 1 - \frac{\sigma^2g^2}{\left(1 - \frac{1}{\delta}\right)^2 - \hat{a}^2},
\end{align}
where $\hat{a} = a - \frac{1}{\delta} - \lambda_{sup}g$, $\hat{a}^2 < \left(1 - \frac{1}{\delta}\right)^2$, and $\lambda_{sup} := \underset{\lambda \in \{\lambda_2,\lambda_N\}}{argmax}\Big|\lambda+\bar{\gamma}\tau-\frac{a_0}{g}\Big|$. Furthermore, $\tau := \frac{\lambda_{N_U}}{\lambda_{N_U}+\lambda_{2_D}}$, where $\lambda_{N_U}$ is the maximum eigenvalue of ${\cal L}_U$ and $\lambda_{2_D}$ is the second smallest eigenvalue of ${\cal L}_D$.
\end{definition}

\subsection{System Theoretic Interpretation of Syncronization Condition and Margin} 
The first step towards the system theoretic interpretation to margin definition is to show the stability condition derived in Theorem \ref{1D_suff_original}, i.e., Eq. (\ref{sync_ric}), has the following equivalent,
\begin{eqnarray}
\left(p - \frac{1}{\delta}\right)\left(\frac{1}{p} - \frac{1}{\delta}\right) > \alpha_0^2 \iff \left(1 - \frac{1}{\delta}\right)^2 > E_{\xi}[\left(\hat{a}-\xi g\right)^2]. \label{equivalence}
\end{eqnarray}
The main point of the equivalence is the equivalent stability condition on the right-hand side is independent of $p$. We observe from \eqref{suff_main1}, if $p=q$ is a valid solution, so is $p=\frac{1}{q}$. Applying the inequality of arithematic and geometric mean (AM-GM Inequality) to $p$ and $\frac{1}{p}$, from Eq. \eqref{suff_main1} we obtain,
\begin{align}\label{suff_main2}
\left(1 - \frac{1}{\delta}\right)^2 = 1 - \frac{2}{\delta} + \frac{1}{\delta^2} \geq 1 - \left(p + \frac{1}{p}\right)\frac{1}{\delta} + \frac{1}{\delta^2} > \alpha_0^2.
\end{align}
We know, if $p = q > 1$ is a solution of \eqref{sync_ineq}, then $p = \frac{1}{q}$ and $p = 1$ are also solutions, since we obtain from \eqref{suff_main2},
\begin{align} 
\left(1 - \frac{1}{\delta}\right)^2 > \left(p - \frac{1}{\delta}\right)\left(\frac{1}{p} - \frac{1}{\delta}\right) > E_{\xi}[\left(\hat{a}-\xi g\right)^2].
\end{align}
Hence, 
we have for some $r > 0$,
\begin{align}\label{suff_condn} 
\left(1 - \frac{1}{\delta}\right)^2 = E_{\xi}[\left(\hat{a}-\xi g\right)^2] + r > E_{\xi}[\left(\hat{a}-\xi g\right)^2] + \frac{r}{2}.
\end{align}
Now, consider some $\epsilon > 0$, such that
\begin{align}\label{eps_val}
\frac{r}{2} = \left(\frac{\epsilon^2}{1+\epsilon}\right)\frac{1}{\delta} = \left(1+\epsilon + \frac{1}{1+\epsilon} - 2\right)\frac{1}{\delta}.
\end{align}
From \eqref{eps_val}, we obtain,
\begin{align}\label{nec_suff}
\left(1+\epsilon - \frac{1}{\delta}\right)\left(\frac{1}{1 + \epsilon} - \frac{1}{\delta}\right) = 1 +\frac{1}{\delta^2} - \frac{2}{\delta} - \left(1+\epsilon+\frac{1}{1+\epsilon}-2\right)\frac{1}{\delta} = \left(1 - \frac{1}{\delta}\right)^2 - \frac{r}{2}.
\end{align}
Using \eqref{suff_condn}, and \eqref{nec_suff}, and setting $p = 1 + \epsilon > 1$, we obtain
\begin{align}
\left(p - \frac{1}{\delta}\right)\left(\frac{1}{p} - \frac{1}{\delta}\right) = \left(1 - \frac{1}{\delta}\right)^2 - \frac{r}{2} > E_{\xi}\left[\left(\hat{a}-\xi g\right)^2\right].
\end{align}
Hence, \eqref{suff_main2} is a necessary and sufficient condition for \eqref{sync_ric}, implying equivalence (\ref{equivalence}).

\begin{figure}[ht!]
\begin{center}
\includegraphics[width = 4.5in]{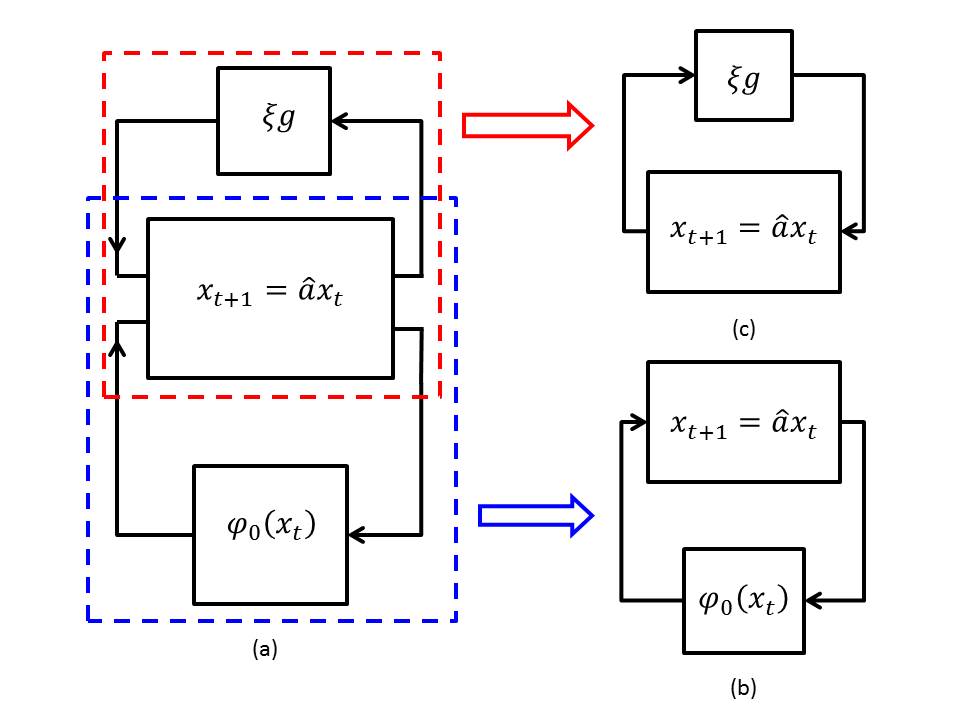}
\caption{(a) System with $H_{\infty}$-norm bounded nonlinearity and stochastic uncertainty in feedback, (b) system with $H_{\infty}$-norm bounded nonlinearity in feedback, and (c) system with stochastic uncertainty in feedback.}
\label{blk_diag}
\end{center}
\end{figure}

One observes the sufficient condition, as provided in Theorem \ref{1D_suff_original} (i.e., Eq. \eqref{sync_ric}), is a Riccati equation in one dimension. For the scalar i.i.d. random variable, $\zeta$, writing $E[\zeta] = \mu := \lambda_{sup}$, $E[(\zeta - \mu)^2] = \sigma^2 := 2\bar{\gamma}\tau\lambda_{sup}$, we can write Eq.  \eqref{sync_ric} as
\begin{align}\label{sync_ric_mod}
p &>  E_{\zeta}\left[(a_0 - \zeta g)^2p + (a_0 - \zeta g)^2\frac{p^2}{\delta - p} + \frac{1}{\delta}\right].
\end{align}
In fact, the above condition is a sufficient condition for stability of the following scalar nonlinear system.
\begin{align}\label{sys_net}
x_{t+1} = \left(a_0 - \zeta g\right)x_t + \phi_0\left(x_t\right) = \left(a_0 - \mu g\right)x_t + \phi_0\left(x_t\right) - \xi gx_t,\qquad\quad \frac{1}{\delta} > \|\phi_0(x)\|_{\infty},
\end{align}
where $\|\phi_0(\cdot)\|_{\infty}$ is the $H_{\infty}$ of $\phi_0$ (\cite{Paganini_book}) and $\xi = \zeta - \mu$ is a zero mean random variable with variance $\sigma^2$. The CoD for $\zeta_c$ is $\gamma_c = \frac{\sigma_c^2}{\mu_c} = 2\bar{\gamma}\tau$. Equation \eqref{sync_ric_mod} can be rewritten as
\begin{align}\label{sync_ineq}
\left(p - \frac{1}{\delta}\right)\left(\frac{1}{p} - \frac{1}{\delta}\right) > E_{\zeta_c}\left[\left(a_0 - \zeta_c g\right)^2\right] = (a_0 - \mu_c g)^2 + \sigma_c^2g^2.
\end{align}
We would like to relate the above stability condition to results from robust stability theory. The central premise of this theory is, if the product of the two gains consisting of a system in the forward loop and the feedback loop is less than unity, then the feedback interconnection is stable. The product of the two gains is referred to as loop gain. In Fig. \ref{blk_diag}(a),  we represent the dynamics of a scalar system (\ref{sys_net}) as a feedback interconnection of a mean linear system with two feedback loops consisting of nonlinearity and stochastic uncertainty. In the following, we analyze each of these feedback loops separately (as shown in Figs. \ref{blk_diag}(b) and \ref{blk_diag}(c)) to derive the stability condition in terms of loop gain. Then, we combine the two separate stability conditions to show the main result of this paper can alternatively be interpreted in terms of loop gains, thereby leading to the proposed definition of a synchronization margin. 

{\noindent \bf Robust Stability to Norm-bounded Nonlinearity:} In Fig. \ref{blk_diag}(b), we have a scalar linear system in the forward loop and a norm bounded nonlinearity in the feedback loop. The scalar system is given by
\begin{align}\label{mean_scalar}
G:\quad x_{t+1} = \hat{a}x_t := (a_0 - \mu g)x_t,
\end{align}
where $a_0 = a - \frac{1}{\delta}$ is the dynamics of the linear system with mean uncertainty in feedback and $G$ represents the system with the transfer function, $G(z) = \frac{1}{z-\hat{a}}$. In our context, the feedback loop with gain, $g$, comes from the nominal network, but is not shown in the figure for simplicity of explanation. Furthermore, the linear dynamics is assumed stable. Results from robust stability and Small Gain Theorem \cite{Astrom_book, Paganini_book} can be used to study the stability of the closed loop. If $\|G\|_{\infty}$ represents the $H_{\infty}$ norm of the system and $\|\phi_0\|_{\infty}$ denotes the $H_{\infty}$ norm of the feedback nonlinearity, the closed loop is stable if
\begin{align}\label{robust_stab}
1 > \|G\|_{\infty}\|\phi_0\|_{\infty} = \left(\frac{1}{1 - |\hat{a}|}\right)\frac{1}{\delta}, \qquad {\rm where} \qquad \|G\|_{\infty} = \frac{1}{1 - |\hat{a}|},\;\;\;\;\;\parallel \phi_0 \parallel_\infty=\frac{1}{\delta} .
\end{align}
The condition for robust stability provided in \eqref{robust_stab} can be reformulated as
\begin{align}\label{robust_stab1}
\left(1 - \frac{1}{\delta}\right)^2 > \hat{a}^2.
\end{align}

{\noindent \bf Robust Stability to Stochastic Uncertainty:} Using results developed in \cite{scl04} under a setting of a linear time invariant system with vector states, we can analyze the stability of the feedback interconnection of a scalar linear system with stochastic uncertainty (as shown in Fig. (\ref{blk_diag})c). The condition for mean square stability of the closed loop system, as obtained from the Small Gain Theorem in \cite{scl04}, is
\begin{align}\label{ms_stab}
1 > \|G\|_{MS}\|\xi g\|_{MS},
\end{align}
where $\|G\|_{MS}$ and $\|\xi g \|_{MS}$ denote the mean square norm of the system and the stochastic uncertainty, respectively. It is shown in \cite{scl04}, that $\|G\|_{MS} = \|G\|_2^2$ for a scalar system, where $\|G\|_2$ denotes the $H_2$-norm of the transfer function $G(z)$. The $H_2$-norm of the scalar system $G(z)$ with an impulse response $h(k) = \hat{a}^k$ is given by
\begin{align}
\|G\|_2 = \left(\sum_{k=0}^{\infty}h(k)^2 \right)^{\frac{1}{2}} = \left(\sum_{k=0}^{\infty}\left(\hat{a}^k\right)^2 \right)^{\frac{1}{2}}.
\end{align}
Thus, the mean square norm of the the linear system is given by
\begin{align}\label{ms_norm}
\|G\|_{MS}^2 = \sum_{k=0}^{\infty}\left(\hat{a}^2\right)^k = \frac{1}{1 - \hat{a}^2}.
\end{align}
The mean square norm of the stochastic uncertainty with zero mean is simply given by its variance, i.e., $\|\xi g\|^2_{MS} = \sigma^2 g^2$. The stability condition \eqref{ms_stab} for mean square stability of the feedback loop for the system in Fig. (\ref{blk_diag})c can now be formulated as
\begin{align}\label{ms_stab1}
 \quad 1 > \left(\frac{1}{1 - \hat{a}^2}\right)\sigma^2g^2 \quad \iff \quad 1 > \hat{a}^2 + \sigma^2 g^2 =  E_{\xi}\left[ (\hat{a} - \xi g)^2\right].
\end{align}

{\noindent \bf Robust Stability to Norm-bounded Nonlinearity and Stochastic Uncertainty:} Stability conditions from the two individual feedback loops as discussed in the previous two sections can be combined to obtained the stability condition in terms of loop gain for the entire system as  shown in Fig. \ref{blk_diag}(a). In particular, the condition for the mean square synchronization of the network, where both the norm-bounded nonlinearity and the stochastic uncertainty are present, can be written as 
\begin{align}\label{ms_sync}
\left(1 - \frac{1}{\delta}\right)^2 > \hat{a}^2+\sigma^2 g^2.
\end{align}
The stability condition  in Eq. (\ref{ms_sync}) is the correct method for combining the stability condition for the two individual feedback loops as expressed in Eqs. (\ref{robust_stab1}) and (\ref{ms_stab1}) for the following reasons. Using $E_\xi[(\hat a-\xi g)^2]=\hat a^2+\sigma^2 g^2$ and the equivalent relationship (\ref{equivalence}), we notice (\ref{ms_sync}) is exactly the mean square synchronization condition as derived in the main results of the paper. Furthermore, when $\sigma=0$ (i.e., no stochastic uncertainty), condition (\ref{ms_sync}) reduces to condition (\ref{robust_stab1}). Similarly, when $\delta=\infty$ (i.e., zero nonlinearity), then condition (\ref{ms_sync}) reduces to (\ref{ms_stab1}). Hence, Eq. (\ref{ms_sync}) can be viewed as the proper generalization of the stability conditions from the two individual loops to the entire system with these two loops operating in tandem.

{\noindent \bf Discussion on Synchronization Margin:} Mean square synchronization condition as expressed in Eq. (\ref{ms_sync}) can be used to derive the expression for mean square synchronization margin. In particular, we note condition (\ref{ms_sync}), after algebraic manipulation, can equivalently be written as  
\begin{align}\label{equiv2}
\left(1 - \frac{1}{\delta}\right)^2 > \hat{a}^2+\sigma^2 g^2 
\quad \iff \quad  1 > \left(\frac{1}{\left(1 - \frac{1}{\delta}\right)^2 - \hat{a}^2}\right)\sigma^2g^2. 
\end{align}
The equivalent condition (\ref{equiv2}) has a nice system theoretic-based interpretation in terms of loop gain. In particular, the quantity, $\left(\frac{1}{\left(1 - \frac{1}{\delta}\right)^2 - \hat{a}^2}\right)$, can be thought of as the gain of the mean linear system with nonlinearity in the feedback loop and $g^2 \sigma^2$ as the gain of the stochastic uncertainty. The system will have a larger margin of stability, if the product of these two quantities is further from one. In particular, the smaller the quantity $\left(\frac{1}{\left(1 - \frac{1}{\delta}\right)^2 - \hat{a}^2}\right)$, the greater the variance, $\sigma^2 g^2$, that can be tolerated to maintain stability, and, hence, more robust the system is to stochastic uncertainty and vice versa. This motivates us to propose the following definition of synchronization margin as a quantity, which measures how far the loop gain is from one.
\begin{align}
\rho_{SM}=1-\left(\frac{1}{\left(1 - \frac{1}{\delta}\right)^2 - \hat{a}^2}\right)\sigma^2g^2.
\end{align}
This is precisely the definition of synchronization margin as proposed in Definition \ref{def_sm}.

\subsection{Interplay of Internal Dynamics, Network Topology and Uncertainty Characteristics}

\begin{figure*}[ht!]
\begin{center}
\includegraphics[width = 0.9\textwidth]{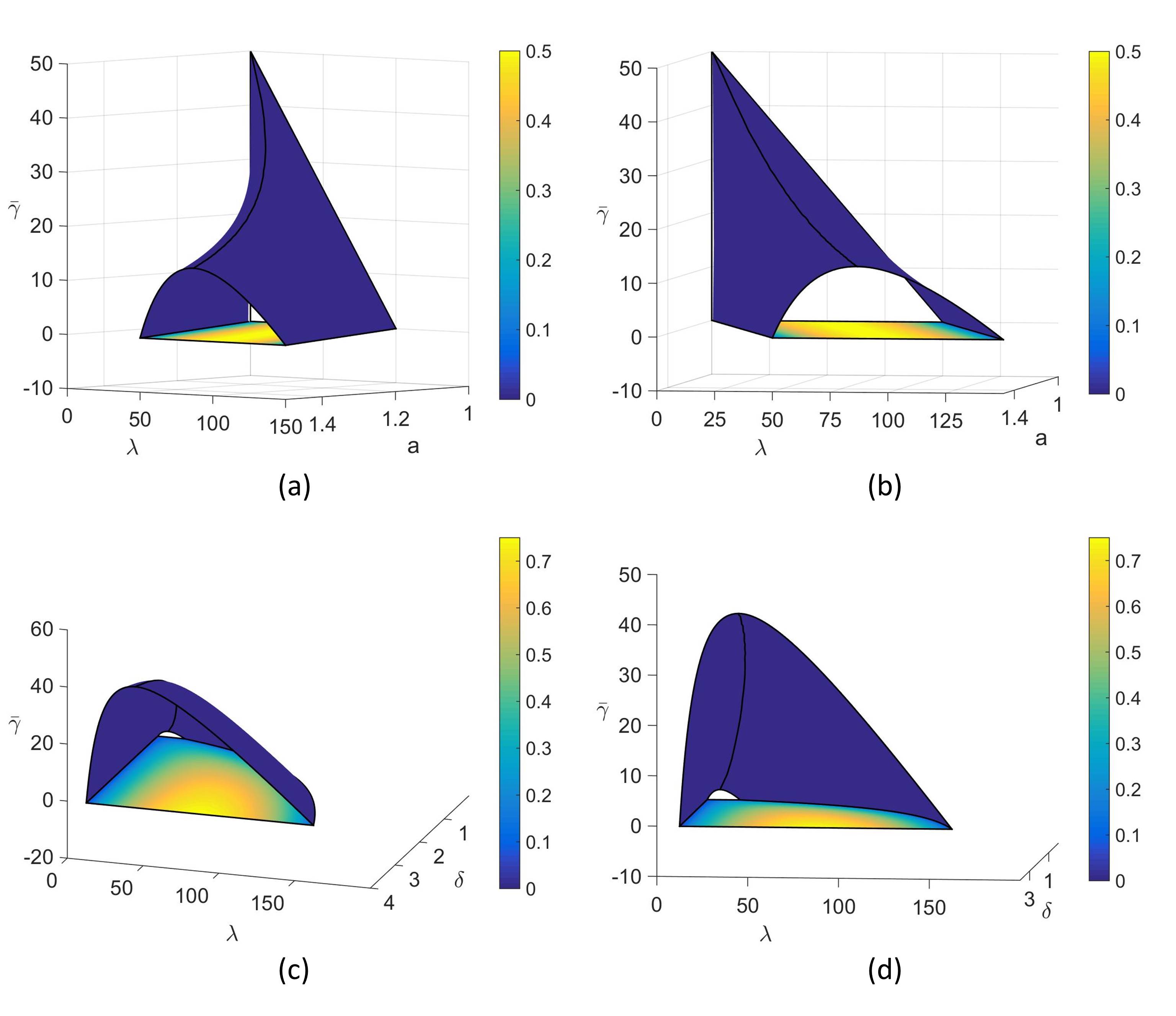}
\caption{(a)\&(b) $a-\lambda-\bar \gamma$ parameter space indicating $\rho_{SM}$ for $g=0.01$, and $\delta=2$, (c)\&(d) $\delta-\lambda-\bar \gamma$ parameter space indicating $\rho_{SM}$ for $a=1.125$ and $g=0.01$.}
\label{fig3d}
\end{center}
\end{figure*}

We now study the interplay of the internal dynamics ($a$), nonlinearity bound ($\delta$), network topology ($\lambda$), and the uncertainty characteristics ($\bar{\gamma}$) through simulations, using a set of parameter values. To nullify the bias of uncertain link locations, we choose to work with a large number of uncertain links to obtain $\tau \approx 1$. In Fig. \ref{fig3d}, we provide different orientations of the 3-dimensional plots used to discuss the interplay between various system, network, and uncertainty parameters over a network with $1000$ nodes. In Figs. \ref{fig3d}(a) and \ref{fig3d}(b), we plot the boundary for the region with a positive $\rho_{SM}$ in the $a-\lambda-\bar{\gamma}$ space, for $g=0.01$, and $\delta = 2$. In Figs. \ref{fig3d}(c) and \ref{fig3d}(d), we plot the boundary for the region with a positive $\rho_{SM}$ in the $a-\lambda-\bar{\gamma}$ space, for $g=0.01$, and $a = 1.125$.

\subsection{Location of Uncertainty and Network Topology}
In the main document, we briefly discussed the impact of the parameter, $\tau$, on the synchronization margin, $\rho_{SM}$. There exists an inverse relation between $\tau$ and $\rho_{SM}$, given by $\rho_{SM} = 1 - \tau\left(\frac{2\lambda_{sup}\bar{\gamma}g^2}{\left(1 - \frac{1}{\delta}\right)^2 - \left(a - \frac{1}{\delta} - \lambda_{sup}g\right)^2}\right)$, thus, making higher $\tau$ detrimental for robustness (low value of $\rho_{SM}$). In this section, we will study the interplay between network topology and the location of uncertainty within the network. In particular, we wish to analyze the average robustness of Small World networks to random network interconnections becoming uncertain. 

\begin{figure*}[ht!]
\begin{center}
\includegraphics[width = 0.8\textwidth]{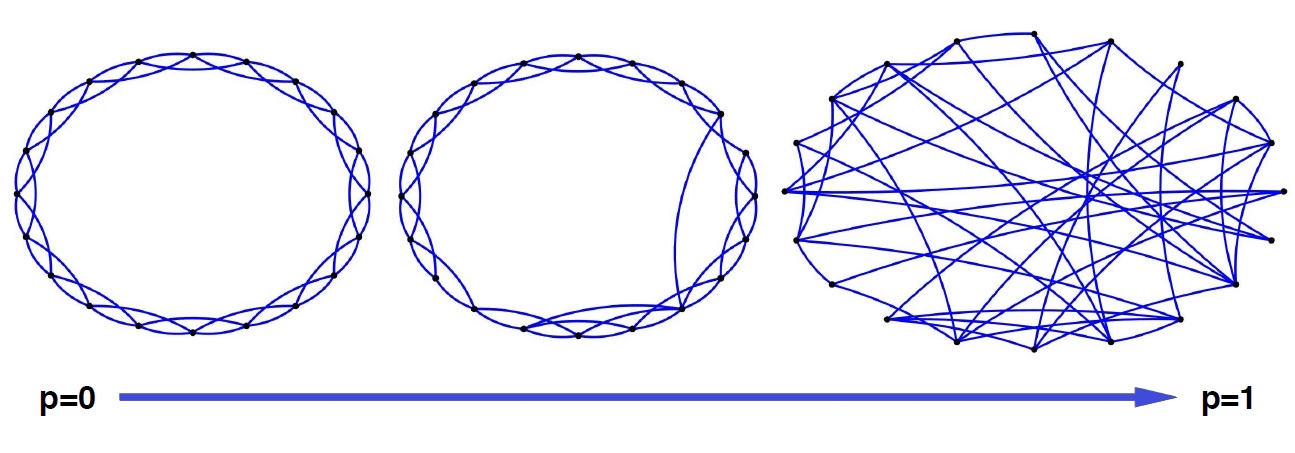}
\caption{Small World network connectivity graph with increasing value of rewiring probability $p$.}
\label{figSM}
\end{center}
\end{figure*}

\begin{figure*}[ht!]
\begin{center}
\includegraphics[width = 0.7\textwidth]{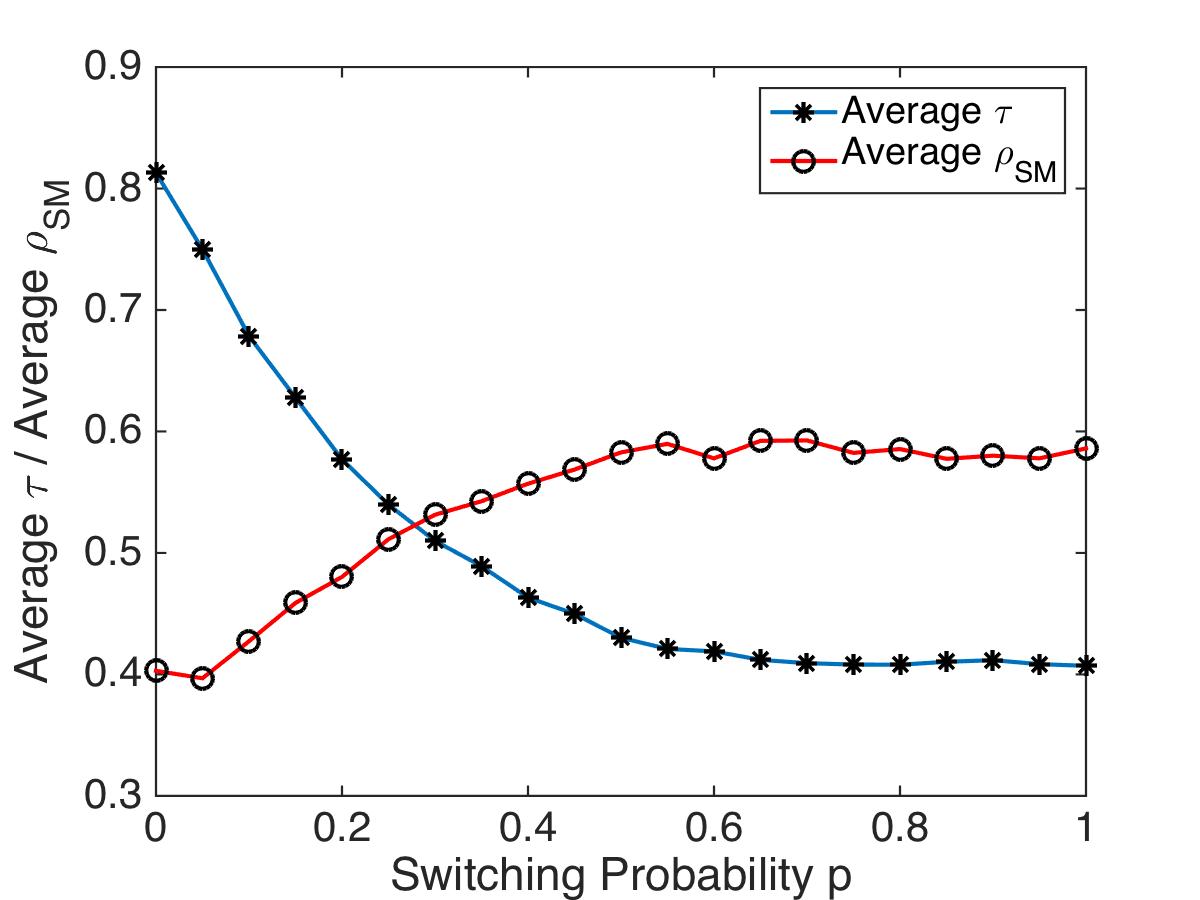}
\caption{Small World network average $\bar{\tau}_{avg}$ as a function of rewiring probability $p$.}
\label{fig_tau}
\end{center}
\end{figure*}

Small World networks were first introduced in \cite{Nature_watts}, and constructed from nearest neighbor networks with random rewiring of links between nodes with a chosen probability. When the rewiring probability, $p=0$, the network is a nearest neighbor network. As $p$ increases, the network loses its nearest neighbor property and has an increasing number of long distance interconnections. This can be visualized through the schematic in Fig. \ref{figSM}, which shows the change in the network connectivity as the rewiring probability increases from $p=0$ to $p=1$ \cite{Nature_watts}.

We initially consider a network of 50 nodes with 8 nearest neighbors per node. Then, we increase the rewiring probability from $p=0$ to $p=1$ in steps of $0.1$ to obtain various random networks. For each such network, we cycle through all the individual links making them uncertain. Then, we compute the value of $\tau$ corresponding to the particular uncertain link with a unit mean value for the interconnection weight. These values of $\tau$ are used to find the average value of $\tau_{avg}$ for the given network. Since these networks are random in nature and the interconnections are formed by probabilistically rewiring links from a nearest neighbor network, we obtain the $\tau_{avg}$ values for 50 samples of random networks for a chosen probability, $p$, which are then used to estimate the mean value of $\tau_{avg}$ given by $\bar{\tau}_{avg}$. Then we plot the $\bar{\tau}_{avg}$ values as a function of the rewiring probability in Fig. \ref{fig_tau}. A curve connecting the data points (blue line in Fig. \ref{fig_tau}) is used to indicate the trend.

We plot the trend of the average synchronization margin, $\rho_{SM}$, in Fig. \ref{fig_tau} (black markers with a red curve connecting them that indicates the trend). It can be observed, as the value of average $\tau$ decreases, the synchronization margin increases, indicating increased randomness (high rewiring probability) in network interconnections makes the network more robust to link uncertainty. It should be noted, there are small variations in the trend shown by $\rho_{SM}$, which we attribute to the computation of the average value for $\tau$, $\lambda_2$, and $\lambda_N$ for 50 samples of Small World networks. Overall, this trend suggests as the rewiring probability in increased, the network robustness increases. Thus, we conclude nearest neighbour and Small World networks are less robust to the injection of uncertainty at a random link location in the network, compared to random network.

\subsection{Optimal Gain Result}

In this subsection we provide the lemma and proof for the results, which provides a method to design the optimal coupling gain for synchronization.

\begin{lemma}\label{opt_gain}
For the network system in Eq. \eqref{coupled_dynamics} with SM given by Eq. \eqref{sync_margin}, the optimal gain, $g^*$, to achieve maximum SM is
\begin{align}
g^* = \frac{2(a-\frac{1}{\delta})}{\max\{\lambda_N,\lambda_2+2\bar{\gamma}\tau \}+\lambda_2 + 2\bar{\gamma}\tau}.
\end{align}
\end{lemma}

\begin{proof}
We observe from \eqref{suff_main2} to maximize the synchronization margin with respect to the coupling gain, $g$, we must minimize $\alpha_0^2$ with respect to $g$, and maximize $\alpha_0^2$ with respect to $\lambda$. This is a regular saddle-point optimization problem \cite{Boyd_book}. Hence, for a given $\lambda$, 
\begin{align}
\frac{\partial \alpha_0^2(\lambda,g)}{\partial g} = -2a_0\lambda + 2\left(\lambda^2 + 2\bar{\gamma}\tau\lambda\right)g = 0.
\end{align}
This provides us with optimal gain, and the corresponding $\alpha_0^2$, \begin{align}
g^*(\lambda) = \frac{a_0}{\lambda + 2\bar{\gamma}\tau}, \qquad \text{and}, \qquad  \alpha_0^2(\lambda,g^*(\lambda)) = \frac{2\bar{\gamma}\tau a_0^2}{\lambda + 2\bar{\gamma}\tau}.
\end{align}
The only important eigenvalues for the graph Laplacian that provide limitations on synchronization (small magnitude of $\rho_{SM}$) are $\lambda_2$ and $\lambda_N$. Hence, we obtain,
\begin{align}
g^*(\lambda_2) = \frac{a_0}{\lambda_2 + 2\bar{\gamma}\tau}, \qquad \text{and}, \qquad g^*(\lambda_N) = \frac{a_0}{\lambda_N + 2\bar{\gamma}\tau}.
\end{align}
Since $\lambda_N \geq \lambda_2$, we have
\begin{align}
g^*(\lambda_2) \geq g^*(\lambda_N) \qquad \text{and} \qquad \alpha_0^2(\lambda_2,g^*(\lambda_2)) \geq \alpha_0^2(\lambda_N,g^*(\lambda_N)).
\end{align}
There also exists a value of gain, $g_e$, which provides the exact same synchronization margin for both $\lambda_2$ and $\lambda_N$. This is obtained by equating, 
\begin{align}
\alpha_0^2(\lambda_2,g_e) = \alpha_0^2(\lambda_N,g_e),
\end{align} 
which provides
\begin{align}
\lambda_2^2 g_e^2 + 2\lambda_2\bar{\gamma}\tau g_e^2 - 2a_0\lambda_2 g_e = \lambda_N^2 g_e^2 + 2\lambda_N\bar{\gamma}\tau g_e^2 - 2a_0\lambda_N g_e.
\end{align}
For $\lambda_N \neq \lambda_2$ and $\bar{\lambda} = \frac{\lambda_2 + \lambda_N}{2}$, this gives
\begin{align}
g_e = \frac{a_0}{\bar{\lambda} + \bar{\gamma}\tau}.
\end{align}
Furthermore, the $\alpha_0^2$ value for $g_e$, is given by 
\begin{align}
\alpha_0^2(\lambda_2,g_e) = \alpha_0^2(\lambda_N,g_e) = a_0^2 - \frac{4\lambda_2\lambda_N a_0^2}{\left(\lambda_N + \lambda_2 + 2\bar{\gamma}\tau\right)^2}.
\end{align}
Since $\lambda_N \geq \lambda_2$, we have $g_e \geq g^*(\lambda_N)$. Furthermore, 
\begin{align}
\alpha_0^2(\lambda_N,g_e) \geq \alpha_0^2(\lambda_N,g^*(\lambda_N)) \qquad \text{and} \qquad \alpha_0^2(\lambda_2,g_e) \geq \alpha_0^2(\lambda_2,g^*(\lambda_2)).
\end{align}
We also conclude, $g^*(\lambda_2) \geq g_e$, iff $\lambda_N \geq \lambda_2 + 2\bar{\gamma}\tau$ and $g_e \geq g^*(\lambda_2)$, iff $\lambda_2 + 2\bar{\gamma}\tau \geq \lambda_N$. We observe, $\lambda_N \geq \lambda_2 + 2\bar{\gamma}\tau$, iff
\begin{align}
\alpha_0^2(\lambda_N,g^*(\lambda_2)) \geq \alpha_0^2(\lambda_N,g_e) \geq \alpha_0^2(\lambda_2,g^*(\lambda_2)).
\end{align}
Hence, $g_e$, being the saddle-point solution, is the optimal gain providing the largest possible $\alpha_0^2(\lambda,g)$ and the smallest $\rho_{SM}$. Similarly, $\lambda_2 + 2\bar{\gamma}\tau \geq \lambda_N$, iff 
\begin{align}
\alpha_0^2(\lambda_2,g_e) \geq \alpha_0^2(\lambda_2,g^*(\lambda_2)) \geq \alpha_0^2(\lambda_N,g^*(\lambda_2)). 
\end{align}
This gives $g^*(\lambda_2)$ as the optimal gain. Furthermore, at the optimal gain, we always have $\lambda_{sup} = \lambda_2$. Defining, $\chi := \max\{\lambda_N,\lambda_2+2\bar{\gamma}\tau\}$, we can write the optimal gain, 
\begin{align}
g^* = \frac{2a_0}{\chi + \lambda_2 + 2\bar{\gamma}\tau}.
\end{align}
Hence, for $\lambda_{sup} = \lambda_2$, we obtain
\begin{align}
\rho_{SM}(g^*) = 1 - \frac{2\bar{\gamma}\tau\lambda_2\left(g^*\right)^2}{\left(1 - \frac{1}{\delta}\right)^2 - \left(a_0 - \lambda_2g^*\right)^2}.
\end{align}
\end{proof}

We will now provide a lemma which will help readers understand the discussion in this paper on the significance of the Laplacian eigenvalues. For $G(V,\Ed)$ with node set $V$ and edge set, $\Ed$, let $\Ed_V$ be all possible connections between nodes in $V$. Then, for $\tilde{\Ed} = \Ed_V \setminus \Ed$, the graph, $\tilde{G} = (V,\tilde{\Ed})$, is the compliment of $G$. Let $0=\lambda_1 < \lambda_2 \leq \cdots \leq \lambda_N$ be the eigenvalues for $G$ and $0=\tilde \lambda_1 < \tilde \lambda_2 \leq \cdots \leq \tilde \lambda_N$ be the eigenvalues for $\tilde G$. We state below, the lemma connecting the eigenvalues of graph, $G$, and its complement, $\tilde{G}$ \cite{Merris1994}.
\begin{lemma}\label{graph_complement}
Let $G \equiv (V,\Ed)$ be a graph on $|V| = N$ nodes. Suppose $\tilde{G} \equiv (V,\tilde{\Ed})$ is the complement of $G$, such that $\tilde{G} = K_N\backslash G$, where $K_N$ is the complete graph on $N$ vertices. Let ${\cal L}_G$ and ${\cal L}_{\tilde{G}}$ be the Laplacian matrices of $G$ and $\tilde{G}$ with eigenvalues, $0=\lambda_1 \leq \lambda_2 \leq \cdots \leq \lambda_N$ and $0=\tilde \lambda_1 \leq \tilde \lambda_2 \leq \cdots \leq \tilde \lambda_N$, respectively. Then, we must have
\begin{align}\label{eigen_condn}
\tilde \lambda_1 = \lambda_1 = 0, \; \tilde \lambda_i = N - \lambda_{N-i+2}, \; \forall \; i\in\{2,\ldots,N\}.
\end{align}
\end{lemma}


\subsection{Simulation Results}\label{section_simulation}

In this subsection, we verify the sufficient condition obtained for mean square synchronization through simulation results. We consider the following 1D system,
\begin{align}
x_{t+1} = ax_t -\phi(x_t) + v_t,\label{it_dynamics}
\end{align}
where $a = 1.125$, $\delta = 8$, and $v_t$ is additive white Gaussian noise with zero mean and variance $\omega^2$. Here, $\phi(x)$ is given by
\begin{eqnarray}
\phi(x_t) =& \frac{\textrm{sgn}(x_t)}{8}\left(s_1\left(|x_t|-\epsilon\right)+\left(s_2^2(|x_t|-\epsilon)^2 + s_3 \right)^{\frac{1}{2}}\right),
\end{eqnarray}
where $s_1=1+m_2$, $s_2 = 1-m_2$, $s_3=4m_2\epsilon^2$, $m_2 = \frac{1}{1+10\epsilon^{0.1}}$, and $\epsilon = 0.3$.  The internal dynamics of the system, as described by Eq. \eqref{it_dynamics}, consists of a  double-well potential, with an unstable equilibrium point at the origin and two stable equilibrium points at $x^* = \pm \epsilon\left(\frac{a-1}{a-2} + \frac{m_2(a-1)}{m_2(a-1)-1}\right) = \pm 0.5237$. So, with no network coupling, i.e., $g=0$, the internal dynamics of the agents will converge to the positive equilibrium point, $x^* > 0$, for positive initial conditions. Similarly, if the initial condition is negative, the systems converge to the negative equilibrium point, $x^* < 0$. The double-well potential system is a prototypical example for modeling synchronization phenomena occurring in the natural sciences and engineering systems. For example, collective motion in molecular dynamics \cite{PNAS_Mezic} and synchronization of generators in the power grid \cite{PNAS_Dorfler} can essentially be modeled using double-well potential.\\

\begin{figure*}[ht]
\begin{center}
\includegraphics[width = 5in]{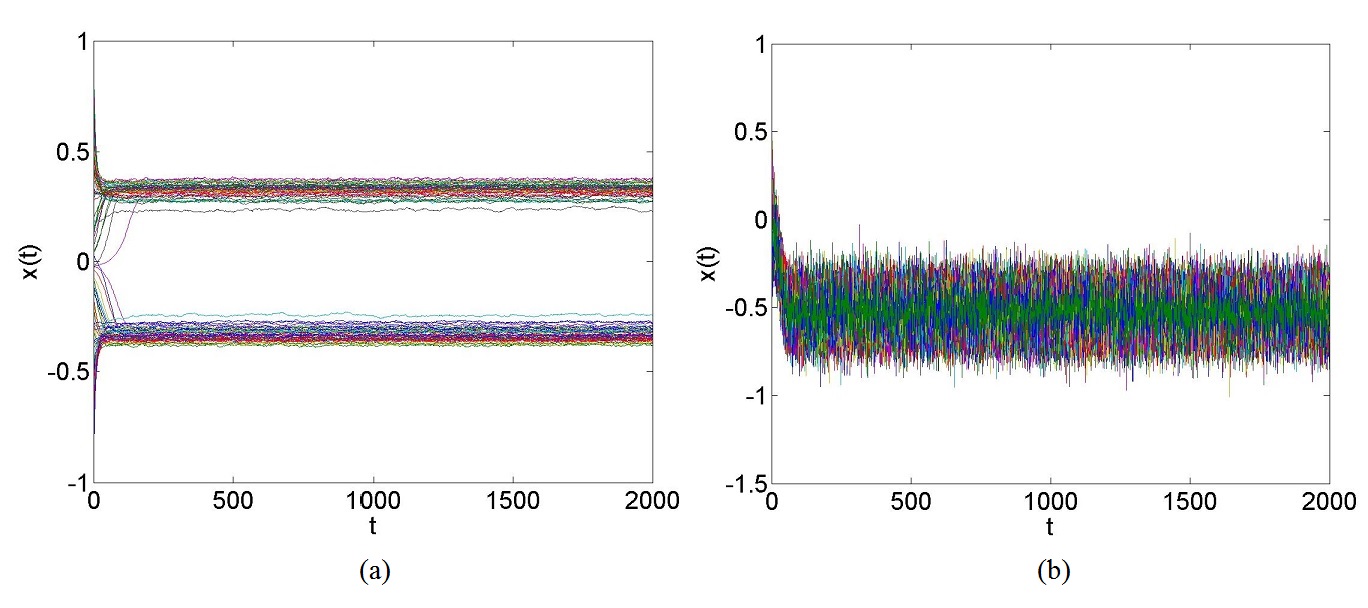}
\caption{(a) Time evolution of systems over a 100-node Small World network, $\bar{\gamma} = 2$, $g = 0.01$, with $\rho_{SM} > 0$, (b) time evolution of systems over a 100-node Small World network, $\bar{\gamma} = 2$, $g=0.001$, with $\rho_{SM} = 0$.}
\label{fig1s3}
\end{center}
\end{figure*}

{\noindent \bf Effect of coupling gain:} We couple this system over a network of 100 nodes, generated as a random network with the Small World property. We choose $60\%$ of the links to be uncertain, making $\tau \approx 1$. The coupling gain for this system is $g = 0.005$. The nominal Laplacian for the network is a standard Laplacian with unit weight. Thus, for all links, $e_{ij}$, connecting nodes $i$ and $j$, $\mu_{ij} = 1$. This network has $\lambda_N = 52.55$ and $\lambda_2 = 26.23$. We now choose $50\%$ of the links in the network to have uncertain weights. The uncertainty in the network link weights is chosen as a uniform variable with zero mean and variance, $\sigma^2 = 2$, such that both these eigenvalues satisfy the required condition from the main result. The CoD of the link uncertainty is $\bar{\gamma} = \frac{\sigma^2}{\mu} = 2$. In Fig. \ref{fig1s3}(a), we plot the results for synchronization of these $100$ systems with simulated additive white Gaussian noise with zero mean and variance, $\omega^2 = 0.1$, which show the systems synchronize in an interval around the equilibrium point.

For systems over the network with identical parameters to those in the previous case and identical link noise variance, if the coupling gain is decreased to $g = 0.001$, which does not satisfy the requirement for the main result, we observe the system is unable to synchronize (Fig. \ref{fig1s3}(b)).\\

\begin{figure*}[ht]
\begin{center}
\includegraphics[width = 7in]{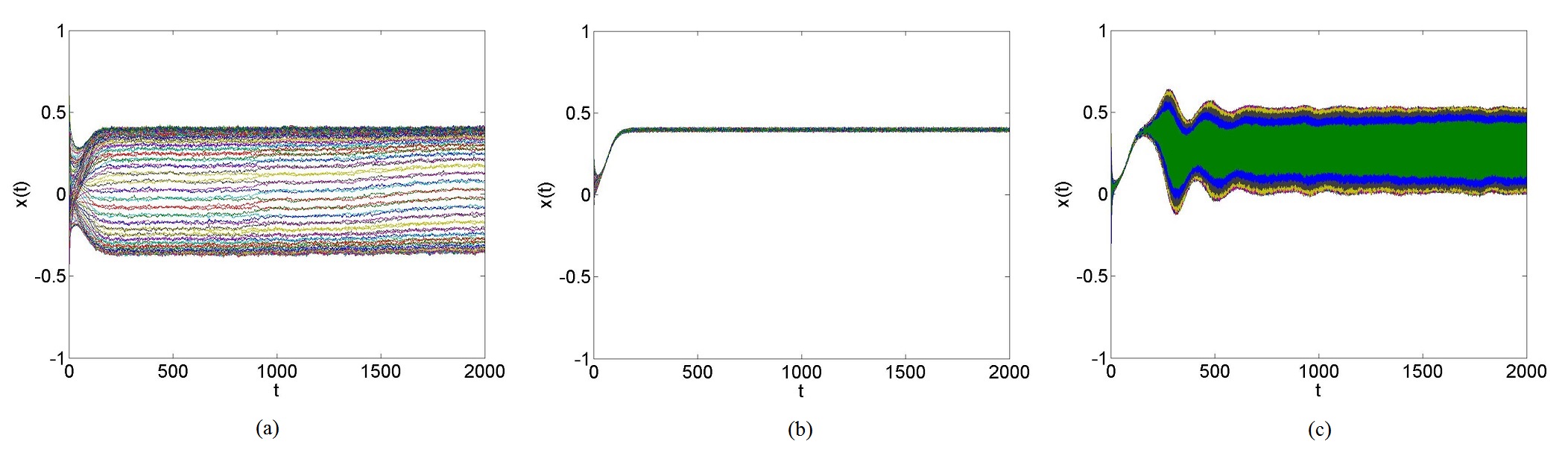}
\caption{(a) Time evolution of systems over a 100-node nearest neighbor network with $6$ neighbors per agent, $\rho_{SM} = 0$, (b) time evolution of systems over a 100-node nearest neighbor network with $20$ neighbors per agent, $\rho_{SM} > 0$, (c) time evolution of systems over a 100-node nearest neighbor network with $32$ neighbors per agent, $\rho_{SM} = 0$.}
\label{fig1s4}
\end{center}
\end{figure*}

{\noindent \bf Effect of number of neighbors:} Next, we study the effect on the group's synchronization ability, due to a change in the number of neighbors for an agent. For this we choose a nearest neighbor network with $100$ nodes. The simulation parameters are chosen as $a=1.05$, $\delta = 16$, $g=0.05$, $\gamma = \frac{1}{12}$, and $\tau \approx 1$. The variance for the uncertainty in the links is chosen small, to clearly observe the effects due to a change in neighbors. As discussed previously, increasing the link uncertainty adds to some numerical inaccuracies in the system causing an additive noise-like effect. Furthermore, the additive noise is also assumed absent to facilitate a clear observation of synchronization.

We first choose a network with $6$ neighbors per agent. For this network we observe the system is unable to synchronize and all the agents break into multiple clusters with each cluster having a small number of agents, Fig. \ref{fig1s4}(a). The agents do not obtain sufficient state information to bind them to the synchronization manifold, due to the small number of neighbors. Now, we increase the number of neighbors to $20$ for each agent. As the number of neighbors increases, the agents synchronize to the synchronization manifold extremely well, with very little noise, Fig. \ref{fig1s4}(b). Furthermore, the rate of synchronization is very high as observed from the simulations, where the agents seem to synchronize within the first $100$ seconds and then collectively move to the synchronization manifold. Synchronization of the agents is observed for a number of neighbors starting at $16$ until the number of neighbors reaches 28. As the number of neighbors increases, we observe significant oscillations before the agents synchronize. Finally, we increase the number of neighbors to $32$ for each agent. This increase in the number of neighbors seems to benefit the synchronization initially, since all agents quickly coalesce together. However, as they approach the synchronization manifold, the high number of neighbors causes the systems to fluctuate significantly about the manifold, leading to an oscillating band of desynchronized agent states, Fig. \ref{fig1s4}(c). This is the beginning of a desynchronized state for the agents. More neighbors for an agent would destabilize the individual system dynamics.